\documentclass[a4paper,12pt]{article}
\usepackage[latin1]{inputenc}

\usepackage[T1]{fontenc}
\usepackage{graphicx}
\usepackage{amsmath,amssymb}
\usepackage{french}
\usepackage{euscript}
\usepackage{enumerate}
\usepackage{epic}
\newcommand{\fp}{\ensuremath{\mathbb{F}}_p}
\newcommand{\fpp}{\ensuremath{\mathbb{F}}_2}
\newcommand{\pp}{\ensuremath{\mathbb{P}}^1}

\newcommand{\dd}{\ensuremath{\mathcal{D}}_0}
\newcommand{\ddd}{\ensuremath{\mathcal{D}}_{0,s}}

\newcommand{\Z}{\ensuremath{\mathbb{Z}}}
\newcommand{\N}{\ensuremath{\mathbb{N}}}

\newcommand{\finetape}{\hfill $\lrcorner$ \medskip}

\newtheorem{theos}{ Théorème  }[section]

\newtheorem{prop}[theos]{Proposition}
\newtheorem{lemme}[theos]{Lemme}
\begin{document}

\title{$\fp$-espaces vectoriels de formes diff\'erentielles
  logarithmiques sur la droite projective}

\author{Guillaume Pagot}

%\affil{Laboratoire de th\'eorie des nombres et algorithmique arithm\'etique,\\ Universit\'e de Bordeaux
%  I, \\ 351, cours de la lib\'eration, \\ 33405 Talence Cedex}

%\email{pagot@math.u-bordeaux.fr}

%\abstract{Soit $p$ un nombre premier, $m$ un entier premier \`a $p$ et $k$ un
%corps alg\'ebriquement clos de caract\'eristique $p$. On \'etudie
%des $\fp$-espaces vectoriels de formes diff\'erentielles logarithmiques ayant chacune
%un seul z\'ero d'ordre $m-1$ en $\infty$ (sauf pour la forme nulle).On
%discute l'existence de tels espaces en fonction de $m$. Nous donnons
%alors une application \`a la d\'eformation en caract\'eristique $0$
%d'action du groupe $(\Z /p\Z)^2$ comme groupe de $k$-automorphismes de
%$k[[t]]$.}

\date{}
\maketitle

\abstract{Let $k$ be an algebraically closed field of characteristic
  $p >0$. Let $m \in \N$, $(m,p)=1$. We study
  $\fp$-vector spaces of logarithmic differential forms on the
  projective line such that each non zero form has a unique zero at
  $\infty$ of given order $m-1$. We
  discuss the existence of such vectors spaces according to the value of
  $m$. We give applications to the lifting to characteristic $0$
  of $(\Z /p\Z)^n$ actions as $k$-automorphisms of $k[[t]]$.}  
%\authorrunninghead{Pagot Guillaume}

%\titleUCfalse

%\authorUCfalse

%\titlerunninghead{$\fp$-espaces vectoriels de formes diff\'erentielles
%  logarithmiques }
{\def\thefootnote{\relax}
\footnote{{\bf Math. Subj. Class. (2000)} : Primary 14D15, 14E22,
  14F10; Secondary 14L30}
\footnote{{\bf Mots cl\'es } : Formes diff\'erentielles logarithmiques
  en caract\'eristique $p$>0, action de $(\Z/p\Z)^n$ sur le disque
  ouvert $p$-adique.}
\begin{center}

\font\twelvebf =cmbx10

\twelvebf{Introduction}

\end{center}

Soit $p$ un nombre premier et $k$ un
corps alg\'ebriquement clos de caract\'eristique $p$. Soit $m$ un
entier premier \`a $p$, l'objet de cet
article est d'\'etudier les $\fp$-espaces vectoriels de dimension
$n\geq 1$ de formes
diff\'erentielles lo\-garithmiques sur $\pp _k$ (i.e. de la forme
$\frac{\mathrm{d}f}{f}$ pour $f\in k(\pp)$), dont les \'el\'ements
non nuls ont un seul z\'ero d'ordre $(m-1)$  en $\infty$. Un tel espace
vectoriel sera not\'e $L_{m+1,n}$. Dans cette \'etude, nous nous
int\'eressons principalement au cas o\`u $n$ est \'egal \`a $2$. Nous
donnons \'egalement des r\'esultats en ce qui concerne les espaces
vectoriels de dimension sup\'erieure.  

Nous commen\c cons par expliciter les conditions impos\'ees aux
formes diff\'erentielles, et nous donnons quelques exemples afin
d'illustrer la richesse et la complexit\'e de tels objets.     

Si $L_{m+1,2}$ est un espace vectoriel comme au-dessus, un lemme
\'el\'ementaire montre que $m+1 \in p \Z$ et donne une premi\`ere
id\'ee sur la r\'epartition des p\^oles des formes diff\'erentielles
non nulles de $L_{m+1,2}$. Nous indiquons une g\'en\'eralisation pour
$n\geq 2$. Comme il est classique d'exprimer \`a partir de l'op\'eration
de Cartier le fait qu'une forme diff\'erentielle soit logarithmique,
nous aboutissons \`a des conditions alg\'ebriques n\'ecessaires et
suffisantes pour l'existence d'espaces $L_{m+1,n}$, qu'il est cependant
difficile d'exploiter.

Nous montrons le th\'eor\`eme suivant qui traite du cas particulier
o\`u $p=2$ et donne une param\'etrisation de tous les espaces $L_{m+1,2}$.

%\newpage

\begin{theos}\label{theo : 0.1}
Supposons $p=2$ et posons $m+1=2n$.

Soit $x_1,\cdots ,x_n \in k$
deux \`a deux distincts, et $u \neq v \in k^*$. Alors il existe
$f_1(z)=\prod_{i=1}^n (z-x_i)(z-y_i)$ et $f_2(z)=\prod_{i=1}^n
(z-x_i)(z-z_i)$ avec
$x_1,\cdots,x_n,y_1,\cdots,y_n,z_1,\cdots,z_n$ deux \`a
deux distincts, tels que les formes diff\'erentielles
$\omega_1:=\frac{\mathrm{d}f_1}{f_1}$ et $\omega_2
:=\frac{\mathrm{d}f_2}{f_2}$ soient de la forme :  

$$\omega_1=\frac{u\mathrm{d}z}{\prod\limits_{i=1}^n(z-x_i)(z-y_i)}
  \;\;\; \mathrm{et}
  \;\;\;\omega_2=\frac{v\mathrm{d}z}{\prod\limits_{i=1}^n(z-x_i)(z-z_i)}$$

Ainsi $\fpp \omega_1+\fpp \omega_2$ est un $L_{m+1,2}$.

R\'eciproquement, tout espace $L_{m+1,2}$ est de cette forme.
\end{theos}

\noindent
Lorsque $p \neq 2$ nous d\'ecrivons les $L_{m+1,2}$ seulement pour $m$ petit.

\begin{theos}\label{theo : 0.2}
On consid\`ere le cas $p \geq 3$. 
\begin{enumerate}[1. ]
\item  Supposons que $m+1=p$. Alors il n'existe pas d'espaces
   vectoriels $L_{m+1,2}$.

 \item Supposons que $m+1=2p$. Alors il existe un espace vectoriel
   $L_{m+1,2}$ si et
   seulement si $p=3$.

 \item Supposons que $m+1=3p$. Alors il n'existe pas d'espaces
   vectoriels $L_{m+1,2}$.
\end{enumerate}
\end{theos}

La d\'emonstration de ce th\'eor\`eme ne fait pas appel \`a
l'op\'eration de Cartier mais \`a une analyse alg\'ebrique des
\'equations sur les r\'esidus aux p\^oles, ce qui la rend technique.
La conclusion d\'epend d'un lemme (lemme \ref{lemme : 2.10})
d'alg\`ebre \'el\'ementaire dont nous n'avons pas vu trace dans la
litt\'erature.
 
Nous donnons \'egalement des exemples de tels espaces
vectoriels de dimension quelconque en suivant une construction d\^ue
\`a Matignon (cf \cite{mat}).

Le th\'eor\`eme \ref{theo : 0.2} a des applications dans le rel\`evement \`a la
caract\'eristi\-que $0$ d'action de groupes. En effet, consid\'erons
un automorphisme $\sigma$ d'ordre $p$ du disque ouvert $p$-adique. On
lui associe naturellement le mod\`ele minimal semi-stable qui
d\'eploie les points fixes de $\sigma$. La fibre sp\'eciale de ce
mod\`ele est un arbre de droites projectives et des formes
diff\'erentielles logarithmiques apparaissent sur les composantes
terminales de cet arbre. Lorsque l'on \'etudie des actions de
$(\Z/p\Z)^2$ sur le disque ouvert $p$-adique, ce sont alors des
espaces vectoriels de telles formes diff\'erentielles qui peuvent
apparaitre. L'application principale est le th\'eor\`eme suivant, qui
donne de nouvelles obstructions au rel\`evement d'actions de groupe et
justifie ainsi l'introduction des espaces $L_{m+1,n}$.

\begin{theos}\label{theo : 0.3}
Soit $\mathrm{G}=(\Z/p\Z)^2$, $p \geq 3$ et $\mathrm{R}$ un anneau de valuation
discr\`ete dominant l'anneau des vecteurs de Witt de $k$. Supposons
que G est un groupe de $k$-automorphismes de $k[[z]]$ et que chacune
des sous-extensions d'ordre $p$ de $k[[z]]^{\mathrm{G}}$ a un
conducteur \'egal \`a $p$ (i.e $\forall \sigma \in G-\{\mathrm{Id}\},\;
v_z(\sigma(z)-z)=p$). Alors, on ne peut pas relever G en un groupe de
R-automorphisme de R$[[Z]]$. 
\end{theos}

Le second inter\^et des espaces $L_{m+1,n}$ r\'eside dans le
th\'eor\`eme suivant :

\begin{theos}\label{theo : 0.4}
On consid\`ere un $L_{m+1,n}$ et une base $\omega_1,\cdots,\omega_n$
de cet espace, chaque $\omega_i$ s'\'ecrivant
$\frac{\mathrm{d}f_i}{f_i}$. Soit $\zeta$ une racine primitive $p$-i\`eme de
l'unit\'e et R=W$(k)[\pi]$ o\`u $\pi ^m:=\lambda:=\zeta -1$, on note K=Frac(R). Alors on peut
trouver $F_i \in $R$[X]$ relevant $f_i$ tels que le produit fibr\'e
des rev\^etements de $\pp_{\mathrm{K}}$ donn\'es par les \'equations
$Y_i^p=F_i(X)$ induisent apr\`es normalisation un rev\^etement de
$\pp_{\mathrm{K}}$ galoisien de groupe $(\Z/p\Z)^n$ ayant bonne
r\'eduction relativement \`a la valuation de Gauss $T:=\pi^{-p}X$. La
fibre sp\'eciale du mod\`ele lisse correspondant est un rev\^etement \'etale, galoisien de
groupe  $(\Z/p\Z)^n$, de la droite affine $\mathbb{A}_k^1$.
\end{theos}  
%Enfin, la derni\`ere partie donne d'autres exemples de tels espaces
%vectoriels de dimension quelconque qui fournissent des r\'ealisations
%de $(\Z/p\Z)^n$ comme groupes d'automorphismes du disque ouvert $p$-adique.

Je souhaite remercier chaleureusement Michel Matignon pour ses
pr\'ecieuses indications et pour sa disponibilit\'e tout au long de
l'avancement de ce travail. 

\section{Pr\'esentation et approche du probl\`eme}

Soit $p$ un nombre premier, $m$ un entier strictement positif, et $k$
un corps alg\'ebriquement clos de caract\'eristique $p$. On fixe une
fois pour toutes un point $\infty$ de la droite projective $\pp_k$.

\vspace{1em}
\noindent
\underline{D\'efinition} : On note $L_{m+1,n}$ un $\fp$-espace
vectoriel de dimension $n\geq 1$ de formes
diff\'erentielles lo\-garithmiques sur $\pp _k$, dont les \'el\'ements
non nuls ont un seul z\'ero d'ordre $(m-1)$  en $\infty$.

\subsection{Espaces $L_{m+1,1}$}

Nous allons exhiber quelques exemples d'espaces $E_{m+1,1}$ ( il est en
effet l\'egitime de s'interroger sur l'existence de tels objets avant
de consid\'erer des espaces de dimension sup\'erieure).

Soit $\fp \omega$ un espace $L_{m+1,1}$ et $z$ un param\`etre de $\pp
_k-\{\infty\}$ tel que $z=0$ n'est pas p\^ole de $\omega$. Ainsi
%Soit $z$ un param\`etre de $\pp_k$. Soit $\omega$ une telle forme
%diff\'erentielle; elle a $m+1$ p\^oles simples, ainsi : 
\font\twelvebf =cmbx10

$$ \omega = \frac{\mathrm{d}f}{f}$$
avec $$ f=\prod_{i=1}^{m+1}(z-x_i)^{h_i} $$ et $ x_i \in
k^*$, $x_i\neq x_j$,  $h_i \in \Z-p\Z$, $ \sum_{i=1}^{m+1}
h_i=0$ mod $p$. Le $(m+1)$-uplet $(h_i)_i$ est appel\'e une
{\twelvebf donn\'ee d'Hurwitz}. Remarquons que $f$ est d\'efinie \`a une
multiplication pr\`es par une puissance $p$-i\`eme et que $h_i
\equiv$ res$_{x_i}\omega$ mod \nolinebreak $p$.

De plus, $\omega$ a un seul z\'ero d'ordre $m-1$ en $\infty$, donc
$\exists u \in k^*$ tel que :
$$\omega=\sum_{i=1}^{m+1}\frac{h_i}{z-x_i}\mathrm{d}z=\frac{u}{\prod\limits_{i=1}^{m+1}(z-x_i)}\mathrm{d}z$$
%Si on exprime la forme $\omega$ en fonction du nouveau param\`etre
%$x:=\frac{1}{z}$, cela s'\'ecrit
%$$\omega
%=\sum_{i=1}^{m+1}\frac{h_ix_i}{1+x_ix}\mathrm{d}x=\frac{ux^{m-1}}{\prod\limits_{i=1}^{m+1}(1+x_ix)}\mathrm{d}x$$

Remarquons que les conditions impos\'ees sur $\omega$ entra\^inent
que $m \notin p\Z$. En effet, supposons que $m \in p\Z$. Vu que
deg($f$)$\in p\Z$, on aurait alors que deg($f'$)$\equiv -1$ mod $p$,
ce qui est impossible.

Si on exprime la forme $\omega$ en fonction du nouveau param\`etre
$x:=\frac{1}{z}$ propice au d\'eveloppement formel, on obtient :
$$\omega
=\sum_{i=1}^{m+1}\frac{-h_ix_i}{1-x_ix}\mathrm{d}x=\frac{-ux^{m-1}}{\prod\limits_{i=1}^{m+1}(1-x_ix)}\mathrm{d}x.$$

Ainsi l'existence d'un $L_{m+1,1}$ est \'equivalente \`a l'existence
d'une solution du syst\`eme :
\begin{gather}
\left\{\begin{array}{ccc}
\sum\limits_{i=1}^{m+1}h_ix_i^{\ell}=0 & \mathrm{pour} & 1\leq \ell \leq m-1 \\
\prod\limits_{i<j}(x_i-x_j) \neq 0 \\
x_i \in k,\; h_i\in \Z-p\Z & &
\end{array}  \tag{*} \right.
\end{gather}

\noindent
\underline{\textit{Remarque}} : Si on fixe les $h_i\in \Z-p\Z$, et si on voit ce
syst\`eme comme un syst\`eme en
les inconnues $x_i$, alors ce syst\`eme est invariant par homoth\'etie et
translation. Cette remarque est essentielle; dans la preuve du
th\'eor\`eme \ref{theo : 0.2}, on sera amen\'e \`a plusieurs reprises \`a
effectuer une translation ``ad\'equate'' sur les $x_i$.

Par la suite, nous serons amen\'es \`a regarder le cas o\`u $m+1 \in p
\Z$. Exa\-minons donc le premier cas $m+1=p$ ($p>2$).
Si on fixe $x_0$ et $x_1$, alors les
\'equations traduisent le fait que le point $(x_2,\cdots,x_m)$
appartient \`a une sous-vari\'et\'e ferm\'ee de
$\mathbb{A}_{\fp}^{m-1}-V(\Delta)$  de dimension $0$ (avec
$\Delta=\prod_{2\leq i<j}(x_i-x_j$), cf.\cite{gre2}). Dans le cas o\`u
une telle vari\'et\'e est non vide on dit que les $h_i$ sont une
donn\'ee d'Hurwitz. Dans \cite{hen} Prop 3.18, Henrio donne
un crit\`ere suffisant sur les $h_i$ pour \^etre une donn\'ee
d'Hurwitz. Malheureusement, dans le cas $m+1=p$ (et plus
g\'en\'eralement dans le cas $m+1 \in p \Z$), ce crit\`ere ne fournit
que des $(m+1)$-uplets $(h_i)_i$ o\`u tous les $h_i$ sont
\'egaux ($h_i=1, \forall i$). N\'eanmoins on peut exhiber d'autres exemples de donn\'ees
d'Hurwitz gr\^ace \`a la remarque suivante : 

On \'ecrit $p-1=d_1d_2$ comme produit de deux entiers sup\'erieurs ou
\'egaux \`a
deux ( il convient de choisir $p >3$ pour que cela soit
possible). Supposons que l'on connaisse une donn\'ee d'Hurwitz
$(h_i)_{0\leq i \leq d_1}$ (donn\'ee par exemple par le crit\`ere
d'Henrio). On a alors un polyn\^ome $f$ de la forme :
$$f=\prod_{i=0}^{d_1}(z-x_i)^{h_i}$$
et tel que
$\omega:=\frac{\mathrm{d}f}{f}=\frac{u\mathrm{d}z}{\prod_{i=0}^{d_1}(z-x_i)}$.

Apr\`es translation \'eventuelle, on peut supposer que $x_0=0$ et
donc :
$$\omega=\frac{u\mathrm{d}z}{zP(z)}\;\;\;\mathrm{avec}\;\;\;P(z):=\prod_{i=1}^{d_1}(z-x_i)$$ 
L'id\'ee est alors de faire un changement de variables $z:=Q(t)$ tel
que $Q'(t)$ divise $f(Q(t))$. Nous allons donner deux exemples de tels
changement de variables et pr\'eciser dans chaque cas les donn\'ees
d'Hurwitz obtenues.

\vspace{1em}
\noindent
\underline{\textit{Exemple 1}} : Prenons le changement de variables
$z:=t^{d_2}$. On obtient alors la forme diff\'erentielle logarithmique :
$$\frac{d_2u\mathrm{d}t}{tP(t^{d_2})}$$
La d\'etermination des donn\'ees d'Hurwitz correspondantes est
fournie par le calcul des r\'esidus de cette forme
diff\'erentielle. On trouve alors le $p$-uplet :
$$d_2h_0,\underbrace{h_1 \cdots h_1}_{d_2 \hbox{\scriptsize~fois}}, \cdots,\underbrace{h_{d_1} \cdots h_{d_1}}_{d_2 \hbox{\scriptsize~fois}}$$

\vspace{1em}

 \noindent
\underline{\textit{Exemple 2}} : Posons cette fois-ci
$z=Q(t)=t^{d_2-1}(t-\alpha)$, o\`u $\alpha$ est choisi tel que
$(t-\frac{d_2-1}{d_2}\alpha )$ divise $z-x_1$ (i.e $Q'(t)$ divise
$f(Q(t))$). Soit $P_1(z)$ et
$P_{\alpha}(t)$ tels que $P(z)=(z-x_1)P_1(z)$ et 
$z-x_1=P_{\alpha}(t)(t-\frac{d_2-1}{d_2}\alpha )^2$. On obtient
alors la forme suivante :
$$\omega=\frac{d_2u\mathrm{d}t}{t(t-\alpha)P_{\alpha}(t)(t-\frac{d_2-1}{d_2}\alpha
  )P_1(t^{d_2-1}(t-\alpha))}$$
Cette fois-ci, la donn\'ee d'Hurwitz prend la forme :
$$h_0,(d_2-1)h_0,\underbrace{h_1 \cdots h_1}_{d_2-2
  \hbox{\scriptsize~fois}},2h_1,\underbrace{h_2 \cdots h_2}_{d_2
  \hbox{\scriptsize~fois}}, \cdots,\underbrace{h_{d_1} \cdots h_{d_1}}_{d_2 \hbox{\scriptsize~fois}}$$

%Ces quelques exemples montrent qu'il existe des formes
%diff\'erentielles ayant les propri\'et\'es susd\'ecrites ( il est en effet
%l\'egitime de s'interroger sur l'existence de tels objets avant d'en
%consid\'erer des espaces vectoriels).

Ces quelques exemples montrent qu'il existe des formes
 diff\'erentielles ayant les propri\'et\'es susd\'ecrites. En fait, des
 calculs men\'es sur ordinateur (pour de petites valeurs de $p$) 
 montrent que beaucoup de $p$-uplets sont des donn\'ees
 d'Hurwitz. La question de d\'eterminer quels sont les $p$-uplets $(h_i)$
 convenables est d\'ej\`a en soi un probl\`eme int\'eressant et
 difficile.

\vspace{1em}
\noindent
\underline{\textit{Remarque}}: Dans ce qui pr\'ec\`ede, on a utilis\'e
soit le param\`etre $z$, soit le param\`etre $x=\frac{1}{z}$ . En fait,
chacune de ces deux \'ecritures a son inter\^et propre. La premi\`ere est
agr\'eable \`a manipuler quand il s'agit de faire un d\'eveloppement
formel et d'exprimer les \'equations en les $x_i$. La seconde est
plus appropri\'ee pour des changements de variables, voire des calculs
de r\'esidus. Par la suite, il nous arrivera de privil\'egier l'une
des deux \'ecritures selon les besoins.

\subsection{Conditions combinatoires pour les $L_{m+1,n}$ ($n\geq 2$)}

En ce qui concerne les espaces $L_{m+1,2}$, on a un lemme combinatoire
qui pr\'ecise l'arrangement des p\^oles des formes diff\'erentielles
non nulles :
\begin{lemme}\label{lemme : 1.5}
  Soit un espace vectoriel $L_{m+1,2}$; alors $m+1 \in p\Z$. De
  plus si on note $ (\omega_1,\omega_2)$ une base de cet espace, alors
  ces deux formes diff\'erentielles ont exactement $\frac{p-1}{p}(m+1)$ p\^oles en commun. 
\end{lemme}

\noindent
\underline{\textit{D\'emonstration}}: Soit $ (\omega_1,\omega_2)$ une base de
l'espace vectoriel en question. On note $(m+1-\lambda)$ le nombre de
p\^oles communs \`a $\omega_1$ et $\omega_2$ (on a donc $0\leq \lambda
\leq m+1$). On note $x_0,\cdots,x_m$ les p\^oles de $\omega_1$, et
$h_0,\cdots,h_m$ les r\'esidus en ces p\^oles. De m\^eme pour
$\omega_2$, on les note
$x_{\lambda},\cdots,x_{\lambda +m}$ et $h_{\lambda}',\cdots,h_{\lambda
  +m}'$ les r\'esidus correspondants (on convient de poser $h_i =0$
pour $i>m$ et $h'_i=0$ pour $i< \lambda$).

Soit $c \in \pp (\fp), \; c=[a,b]$ (en coordonn\'ees homog\`enes);
alors $\omega :=a\omega_1 +b\omega_2$ a exactement $m+1$ p\^oles. Donc,
il existe exactement $\lambda$ valeurs de $i$ pour lesquelles
$ah_i+bh'_i=0$. On a alors partitionn\'e les $(m+1+\lambda)$ points
$x_i$ en $p+1$ ensembles de $\lambda$ points. Ainsi
$(m+1+\lambda)=(p+1)\lambda$ et $m+1=\lambda p$. On v\'erifie
ais\'ement que le nombre de p\^oles communs \`a $\omega_1$ et
$\omega_2$ est celui annonc\'e.

\begin{flushright}
$\square$
\end{flushright}

On peut montrer une g\'en\'eralisation dans le cas des espaces
vectoriels $L_{m+1,n}$ :
\begin{lemme}\label{lemme : 1.6}
On conserve les notations pr\'ec\'edentes. Consid\`erons
un espace vectoriel $L_{m+1,n}$ ($n\geq 2$),
 alors $m+1 \in p^{n-1} \Z$. De plus, si
$(\omega_1,\cdots,\omega_n)$ est une base, alors ces $n$ formes
diff\'erentielles ont exactement $\frac{(p-1)^{n-1}}{p^{n-1}}(m+1)$
p\^oles en commun. 
\end{lemme}

\noindent
\underline{\textit{D\'emonstration}} : La d\'emonstration se fait
par r\'ecurrence sur $n$. On prend donc un $\fp$-espace vectoriel
$L_{m+1,n}$ engendr\'e par $n$ formes diff\'erentielles
lin\'eairement ind\'ependantes
$(\omega_1,\cdots,\omega_n)$. L'hypoth\`ese de r\'ecurrence aux rangs
inf\'erieurs dit que pour $j$ formes diff\'erentielles
($j<n$) parmi les $\omega_i$, ces $j$ formes ont exactement
$\frac{(p-1)^{j-1}}{p^{j-1}}(m+1)$ p\^oles en commun. Notons $T$ le
nombre total des p\^oles apparaissant dans les formes
diff\'erentielles $\omega_1,\cdots,\omega_n$ et
$\lambda$ le nombre de p\^oles communs \`a toutes ces
diff\'erentielles. On note \'egalement
$N_{i_1i_2\cdots i_k}$ le nombre de p\^oles communs aux formes
diff\'erentielles $\omega_{i_1},\cdots,\omega_{i_k}$.Alors, on a la
relation : 
\begin{eqnarray}
T & = & \sum_{k=1}^n(-1)^{k+1}\sum_{i_1<i_2<\cdots<i_k}
N_{i_1i_2\cdots i_k} \nonumber \\
& = & (m+1)\sum_{k=1}^{n-1}(-1)^{k+1}
C_n^k\left(\frac{p-1}{p}\right)^{k-1} +(-1)^{n+1}\lambda \nonumber \\
& = &(m+1)\left(\frac{p}{p-1}\right)\sum_{k=1}^{n-1}(-1)^{k+1}
C_n^k\left(\frac{p-1}{p}\right)^{k} +(-1)^{n+1}\lambda \nonumber
\end{eqnarray}
\begin{eqnarray} 
& =
&-(m+1)\left(\frac{p}{p-1}\right)\left(\left(1-\frac{p-1}{p}\right)^n-1-(-1)^n\left(\frac{p-1}{p}\right)^n\right)
\nonumber \\
 & & +(-1)^{n+1}\lambda
\nonumber \\
& =
&(m+1)\left(\frac{p}{p-1}\right)\left(1+(-1)^{n}\left(\frac{p-1}{p}\right)^n
-\left(\frac{1}{p}\right)^n\right)+(-1)^{n+1}\lambda \nonumber
\end{eqnarray}
On note $x_1,\cdots,x_T$ les p\^oles et $h_{j,i}$ le r\'esidu
(\'eventuellement nul) de la forme diff\'erentielle $\omega_j$ au
point $x_i$. Soit $[a_1,\cdots,a_n] \in \mathbb{P}^{n-1}(\fp)$, alors
on a :
$$a_1h_{1,i}+\cdots +a_nh_{n,i}=0$$
pour exactement ($T-(m+1)$) valeurs de $i$. D'autre part, si on
consid\`ere un point $x_i$, il est p\^ole de toutes les
formes diff\'erentielles sauf celles de la forme
$a_1\omega_1+\cdots+a_n\omega_n$, avec  $a_1h_{1,i}+\cdots
+a_nh_{n,i}=0$ (ce qui fait pour chaque $i$ un total de
($p^{n-2}+p^{n-3}+\cdots +1$) formes diff\'erentielles modulo la multiplication par un
\'el\'ement de $\fp^*$). En r\'esum\'e,
l'ensemble des p\^oles $x_1,\cdots,x_T$ 
est la r\'eunion de ($p^{n-1}+p^{n-2}+\cdots +1$) ensembles de ($T-(m+1)$)
\'el\'ements, chaque \'el\'ement \'etant inclus dans exactement
($p^{n-2}+p^{n-3}+\cdots +1$) de ces ensembles. On a donc la relation
\nolinebreak:
$$T(p^{n-2}+p^{n-3}+\cdots +1)=(T-(m+1))(p^{n-1}+p^{n-2}+\cdots +1)$$
et donc :
$$ T=\frac{(m+1)(p^n-1)}{(p-1)p^{n-1}}$$
En comparant avec l'expression de $T$ d\'ej\`a calcul\'ee
pr\'ec\'edemment, il vient \nolinebreak :
\begin{eqnarray}
\frac{(m+1)p}{p-1}\left[\frac{p^n-1}{p^n}-\left(1+(-1)^{n}\left(\frac{p-1}{p}\right)^n
-\left(\frac{1}{p}\right)^n\right)\right]-(-1)^{n+1}\lambda & = & 0
\nonumber \\
\frac{(m+1)p}{p-1}\left[(-1)^{n+1}\left(\frac{p-1}{p}\right)^n\right]-(-1)^{n+1}\lambda & = & 0
\nonumber \\
\frac{(m+1)(p-1)^{n-1}}{p^{n-1}} - \lambda & = & 0 \nonumber
\end{eqnarray}
Finalement $m+1 \in p^{n-1}\Z$ et $\lambda=\frac{(p-1)^{n-1}}{p^{n-1}}(m+1)$.

\begin{flushright}
$\square$
\end{flushright}

\subsection{Conditions alg\'ebriques sur $L_{m+1,n}$}
%On consid\`ere toujours ($\omega_1$, $\omega_2$) une base d'un espace
%vectoriel $L_{m+1,2}$ et on note $z$ un param\`etre de $\pp_k$. Nous allons montrer que les conditions
%impos\'ees sur $L_{m+1,2}$ peuvent se traduire comme suit :

Soit $z$ un param\`etre de $\pp_k-\{\infty\}$. Nous allons montrer la
proposition suivante :

\begin{prop}\label{prop : 1.7}
Soit $\omega_1$, $\omega_2$ deux formes diff\'erentielles sur $\pp
_k$. Alors $\fp \omega_1+\fp \omega_2$ est un $L_{m+1,2}$ si et
seulement si il existe deux polyn\^omes $A$ et $B$ avec $$\mathrm{deg}(iA+jB)=
\frac{m+1}{p}, \;\;\;\;\;\;\;\forall [i,j] \in \pp(\fp),$$
tels que :
$$\omega_1=\frac{A\mathrm{d}z}{A^pB-AB^p}\;\;\;\mathrm{et}\;\;\; \omega_2=\frac{B\mathrm{d}z}{A^pB-AB^p}$$
et $((A^p-AB^{p-1})^{p-1})^{(p-1)}=-1$.
\end{prop}

\underline{\textit{D\'emonstration :}} Supposons que $\fp \omega_1+\fp \omega_2$ est un
$L_{m+1,2}$. On sait d'apr\`es le lemme \ref{lemme : 1.5} que $m+1=\lambda p$ et
que l'ensemble des p\^oles est partitionn\'e en $p+1$ ensembles de
$\lambda$ p\^oles. Plus pr\'ecis\'ement, on \'ecrit que :
\begin{enumerate}[$\bullet$]
   \item $\omega_1$ a ses p\^oles en les points $x_0,\cdots,x_{\lambda p-1}$

   \item $\omega_2$ a ses p\^oles en les points $x_{\lambda},\cdots,x_{\lambda(p+1) -1}$

   \item quitte \`a renum\'eroter, on peut supposer que
    $\omega_1+i\omega_2$ (pour $i$ variant de $1$ \`a $p-1$) a des
    p\^oles en tous les $x_j$ sauf pour $\lambda i \leq j \leq
    \lambda(i+1)-1$.
\end{enumerate}
 
On note $P_j(z)=\prod_{k=\lambda j}^{ \lambda(j+1)-1} (z-x_k)$. Alors
$\omega_1$ et $\omega_2$ s'\'ecrivent :
$$ \omega_1= \frac{uP_0(z)\mathrm{d}z}{\prod_{j=0}^p
  P_j(z)}\;,\;\;\omega_2 =\frac{vP_p(z)\mathrm{d}z}{\prod_{j=0}^p
  P_j(z)} $$
o\`u $u$ et $v$ sont des constantes non nulles.

\noindent
On a alors deux \'ecritures pour $\omega_1 +i\omega_2$ :
$$\omega_1 +i\omega_2=\frac{(uP_0(z)+ivP_p(z))\mathrm{d}z}{\prod_{j=0}^p
  P_j(x)}= \frac{(w_iP_i(z))\mathrm{d}z}{\prod_{j=0}^p
  P_j(z)}$$
o\`u $w_i$ est une constante non nulle.

\noindent
On a donc $uP_0(z)+ivP_p(z)=w_iP_i(z)$. En identifiant les termes
dominants de chaque expression, on trouve $w_i=u+iv$ et donc
$uP_0(z)+ivP_p(z)=(u+iv)P_i(z)$.

 Le rapport $\frac{u}{v}$ n'est pas dans $\fp$. En effet, si $\frac{u}{v}=-i \in \fp$,
alors $(u+iv)P_i=0=u(P_0-P_p)$ et  donc $P_0=P_p$, ce qui implique que les
$x_j$ ne sont pas distincts.

Posons $a=\frac{u}{v}$. Alors :
$$\omega_2=v\prod_{j=0}^{p-1}\frac{(a+j)}{(aP_0+jP_p)}\mathrm{d}z=\frac{v(a^p-a)}{(aP_0)^p-aP_0P_p^{p-1}}\mathrm{d}z$$
et 
$$\omega_1=a\omega_2\frac{P_0}{P_p}=\frac{v(a^p-a)}{(aP_0)^{p-1}P_p-P_p^{p}}\mathrm{d}z.$$

Soit $\alpha \in k$ tel que $\alpha ^pv(a^p-a)=1$ et posons
$A:=\alpha aP_0$, $B:=\alpha P_p$. Vu que $a\notin \fp$, on a :

$$\mathrm{deg}(iA+jB)=
\frac{m+1}{p}, \;\;\;\;\;\;\;\forall [i,j] \in \pp(\fp),$$

Il reste maintenant \`a exprimer
le fait que les formes diff\'erentielles :
$$\frac{A\mathrm{d}z}{A^pB-AB^p}\;\;\;\mathrm{et}\;\;\; \frac{B\mathrm{d}z}{A^pB-AB^p}$$
sont logarithmiques.
Pour exprimer cette condition, on peut exprimer la relation
$\mathcal{C}\omega_i =\omega_i$, o\`u la lettre $\mathcal{C}$
d\'esigne l'op\'eration de Cartier. Rappelons de quoi il s'agit; si on
consid\`ere une forme diff\'erentielle $\omega$, alors on peut
l'\'ecrire :
$$\omega =(f_0^p(z)+zf_1^p(z)+\cdots+z^{p-1}f_{p-1}^p(z))\mathrm{d}z$$
On d\'efinit $\mathcal{C}\omega=f_{p-1}\mathrm{d}z$. Une
condition n\'ecessaire et suffisante pour que $\omega$ soit
logarithmique est que $\mathcal{C} \omega =\omega$ (dans le cas de
formes diff\'erentielles sur $\pp$, la preuve est \'el\'ementaire). Remarquons que
cette condition de Cartier peut \'egalement s'exprimer de la fa\c con
suivante : si on \nolinebreak a $\omega =f$d$z$ alors $\omega$ est logarithmique si
et seulement si :
$$f^{(p-1)}=-f^p$$
 A l'aide de cette op\'eration, on va
montrer que les hypoth\`eses ``$\omega_1$ est logarithmique'' et
``$\omega_2$ est logarithmique'' sont \'equivalentes. 

Supposons en effet que $\frac{B\mathrm{d}z}{A^pB-AB^p}$ est
logarithmique. En \'ecrivant que :
$$\frac{B\mathrm{d}z}{A^pB-AB^p}=\frac{B(A^pB-AB^p)^{p-1}\mathrm{d}z}{(A^pB-AB^p)^p}$$
on voit que la condition donn\'ee par l'op\'eration de Cartier
s'exprime par l'\'egalit\'e :
$$(B(A^pB-AB^p)^{p-1})^{(p-1)}=-B^p$$
A partir de cette expression, on en tire :
\begin{eqnarray}
(A^pB(A^pB-AB^p)^{p-1})^{(p-1)} & = & -A^pB^p \nonumber \\
(((A^pB-AB^p)+AB^p)(A^pB-AB^p)^{p-1})^{(p-1)} & = & -A^pB^p \nonumber
\\(AB^p(A^pB-AB^p)^{p-1}+(A^pB-AB^p)^p))^{(p-1)} & = & -A^pB^p \nonumber
\\ (AB^p(A^pB-AB^p)^{p-1})^{(p-1)} & = & -A^pB^p \nonumber \\
(A(A^pB-AB^p)^{p-1})^{(p-1)} & = & -A^p \nonumber
\end{eqnarray}
et la derni\`ere \'egalit\'e entra\^ine que
$\frac{A\mathrm{d}z}{A^pB-AB^p}$  est logarithmique.

On peut donc r\'esumer ces conditions en disant que :
%$$\frac{A\mathrm{d}z}{A^pB-AB^p}$$
%est logarithmique.
\begin{gather}
((A^p-AB^{p-1})^{p-1})^{(p-1)}=-1 \tag{**}
\end{gather}

%\finetape

Inversement si on a :
$$\omega_1=\frac{A\mathrm{d}z}{A^pB-AB^p}\;\;\;\mathrm{et}\;\;\;
\omega_2=\frac{B\mathrm{d}z}{A^pB-AB^p}$$
avec $A$,$B$ v\'erifiant les conditions de la proposition, on
montre facilement que les formes $i\omega_1 +j\omega_2$ (pour
$(i,j)\neq 0$) sont logarithmiques et n'ont qu'un seul z\'ero d'ordre
$(m-1)$ en $\infty$.

\begin{flushright}
$\square$
\end{flushright}
 
\vspace{1em}
\noindent
\underline{\textit{Remarque 1 :}}
L'\'equation diff\'erentielle (**) est difficile \`a
manipuler. En effet, si on la d\'eveloppe, il apparait des
d\'eriv\'ees $k$-i\`emes de
puissances de $A$, $A$ \'etant lui-m\^eme de degr\'e $\lambda$ (la
r\'esolution n'apparait simple que dans le cas o\`u $\lambda =1$ ou
$p=2$). 

On peut donner une autre formulation de la condition (**) en termes de
congruence : puisque $f:=A^p-AB^{p-1} \in k[z]$, $\omega_1$ est
logarithmique si et seulement si $(f')^{p-1}=1$ modulo $f$. 

\vspace{1em}
\noindent
%\begin{cors}
% Si $A$ et $B$ sont deux polyn\^omes tels que
% $\frac{\mathrm{d}z}{A^p-AB^{p-1}} $ soit de la forme $
% \frac{\mathrm{d}f}{f}$) (i.e, on a un espace $L_{m+1,2}$), il est alors possible d'en
%construire d'autres gr\^ace \`a certains changements de variables. En
%effet, si on pose $z:=\phi(t):=\alpha t+P(t^p)$ avec $\alpha \in k^*$ et $P \in k[t]$, on trouve
%la forme diff\'erentielle :
%$$\frac{\alpha \mathrm{d}t}{A(\alpha t+P(t^p))^p-A(\alpha t+P(t^p))B(\alpha t+P(t^p))^{p-1}}$$ 
%qui est encore logarithmique puisque elle est \'egale \`a
%$\frac{\mathrm{d}(f\circ \phi)(t)}{(f\circ \phi)(t)}$. Ainsi on
%obtient un espace $L_{(m+1)deg(P)p,2}$.
%\end{cors}

\noindent
\underline{\textit{Remarque 2}} : On a une formulation similaire du
probl\`eme pour les $L_{m+1,n}$
($n\geq 3$). Pour cela, on reprend les notations du lemme \ref{lemme : 1.6}. On
note $P$ un polyn\^ome qui n'a que des racines simples qui sont les p\^oles des
formes diff\'erentielles $\omega_i$. Alors chaque forme $\omega_i$
peut s'\'ecrire $\omega_i=\frac{Q_i}{P}$d$z$ o\`u $Q_i$ est un polyn\^ome
avec pour seules racines simples les points $x_i$ o\`u $\omega_i$ n'a
pas de p\^oles. Pour chaque valeur $[a_1,\cdots,a_n]\in
\mathbb{P}^{n-1}(\fp)$, le polyn\^ome $a_1Q_1+\cdots +a_nQ_n$ a
exactement ($T-(m+1)$) racines simples (toujours parmi les p\^oles des
formes diff\'erentielles ), et chaque point $x_i$ est racines
d'exactement ($p^{n-2}+p^{n-3}+ \cdots +1$) de ces polyn\^omes. On
a donc la relation :
$$P^{(p^{n-2}+p^{n-3}+ \cdots +1)}=\gamma \prod_{i=1}^n \prod_{j_{i-1}=0}^{p-1}
\cdots \prod_{j_1=0}^{p-1} (Q_i+j_{i-1}Q_{i-1}+ \cdots +j_1Q_1)$$
o\`u $\gamma$ est une constante.

\noindent
Quitte \`a multiplier $P$ par une constante, on peut supposer
$\gamma=1$. La condition sur les $\omega_i$ pour \^etre logarithmique
s'exprime en disant que les formes :
$$\frac{P^{p( p^{n-3}+p^{n-4}+ \cdots +1)}Q_i}{\prod_{i=1}^n \prod_{j_{i-1}=0}^{p-1}
\cdots \prod_{j_1=0}^{p-1} (Q_i+j_{i-1}Q_{i-1}+ \cdots +j_1Q_1)}$$
sont logarithmiques. On reconnait au d\'enominateur le d\'eterminant
de Moore des polyn\^omes $Q_1 \cdots Q_n$ (cf. \cite{gos}), ce qui g\'en\'eralise la
forme que l'on avait pour $n=2$; en effet, $A^pB-AB^p$ est le d\'eterminant
de Moore de $A$ et $B$.
\vspace{1em}

\noindent
\underline{\textit{Remarque 3}} : On a vu pr\'ec\'edemment que
lorsqu'on disposait d'un $L_{m+1,2}$ engendr\'e par deux formes
$\omega_1$ et $\omega_2$, les coefficients $u$ et $v$ ``associ\'es''
\'etaient lin\'eairement ind\'ependants sur $\fp$. On peut
g\'en\'eraliser ce r\'esultat aux espaces $L_{m+1,n}$ : soit un espace
$L_{m+1,n}$ engendr\'e par les formes diff\'erentielles
$\omega_1,\cdots,\omega_n$. Comme on l'a vu juste au-dessus on peut
\'ecrire $\omega_i=\frac{Q_i}{P}$d$z$; on choisit de prendre $P$
unitaire et on note $u_i$ le terme de plus haut degr\'e de
$Q_i$. Montrons que les $u_i$ sont lin\'eairement ind\'ependants sur
$\fp$.

\noindent
 Soit $a:=(a_1,\cdots,a_n)\in \fp^n-\{0\}$;
d\'efinissons $\omega_a:=a_1\omega_1+\cdots +a_n\omega_n$. Alors :
$$\omega_a=\frac{a_1Q_1+\cdots+a_nQ_n}{P}\mathrm{d}z:=\frac{Q_a}{P}\mathrm{d}z$$
La forme $\omega_a$ doit avoir le m\^eme nombre de p\^oles que les
$\omega_i$, donc le polyn\^ome $Q_a$ a le m\^eme
degr\'e que les $Q_i$. En particulier le coefficient de plus
haut degr\'e de $Q_a$ est non nul. Donc $a_1u_1+\cdots +a_nu_n \neq
0$.  

\vspace{1em}

\noindent
\underline{\textit{Remarque 4}} : Soit $\Phi$ : $\pp_k \longrightarrow
\pp_k$ donn\'ee par $\Phi(t)=\alpha t+P(t^p)$ avec $\alpha \in k^*$ et
$P \in k[t]$ (i.e $\Phi$ est un rev\^etement \'etale de $\pp_k -\{
\infty\}$). Si $F$ est un $L_{m+1,n}$ engendr\'e par les formes
diff\'erentielles $\omega_1, \cdots ,\omega_n$ (avec
$\omega_i=\frac{\mathrm{d}f_i}{f_i}$) alors $\Phi ^*F$ est un
$L_{(m+1)\mathrm{deg}\Phi ,n}$ (o\`u $\Phi ^*F$ d\'esigne le
$\fp$-espace vectoriel engendr\'e par les formes
$\frac{\mathrm{d}(f_i\circ \Phi)}{f_i\circ \Phi}$).

\section{R\'esultats et applications}

A d\'efaut de pouvoir exploiter la condition (**) d\'ecrite ci-dessus, nous
allons explorer les relations alg\'ebriques entre p\^oles et r\'esidus.
\subsection{Un cas particulier : $p=2$}

Le cas $p=2$ appara\^it comme un cas particulier dans la mesure o\`u
toutes les donn\'ees d'Hurwitz sont \'egales \`a $1$.

%\begin{theo}
%Posons $m+1=2n$ ($m$ est forc\'ement de cette forme d'apr\`es le lemme
%$1.2.1$). Soit $x_1,\cdots ,x_n \in k$ distincts deux \`a deux, et
%$u,v \in k^*$ deux \'el\'ements distincts. Alors il existe
%$y_1,\cdots,y_n$ et $z_1,\cdots,z_n$ dans $k$ tels que les points
%$x_i,y_i,z_i$ soient distincts deux \`a deux et v\'erifient les conditions
%suivantes :

%Si on pose $f_1(x):=\prod_{i=1}^n(1+x_ix)(1+y_ix)$ et
%$f_2(x):=\prod_{i=1}^n(1+x_ix)(1+z_ix)$ alors les formes
%diff\'erentielles $\omega_1 :=\frac{\mathrm{d}f_1}{f_1}$ et $\omega_2
%:=\frac{\mathrm{d}f_2}{f_2}$ sont de la forme :
%$$\omega_1=\frac{ux^{2n-2}\mathrm{d}x}{\prod\limits_{i=1}^n(1+x_ix)(1+y_ix)}
%  \;\;\; \mathrm{et}
%  \;\;\;\omega_2=\frac{vx^{2n-2}\mathrm{d}x}{\prod\limits_{i=1}^n(1+x_ix)(1+z_ix)}$$

%En particulier, l'espace vectoriel engendr\'e par $\omega_1$ et
%$\omega_2$ est de la forme voulue.
%\end{theo}
\begin{theos}\label{theo : 2.8}
Supposons $p=2$ et posons $m+1=2n$.

Soit $x_1,\cdots ,x_n \in k$
deux \`a deux distincts, et $u \neq v \in k^*$. Alors il existe
$f_1(z)=\prod_{i=1}^n (z-x_i)(z-y_i)$ et $f_2(z)=\prod_{i=1}^n
(z-x_i)(z-z_i)$ avec
$x_1,\cdots,x_n,y_1,\cdots,y_n,z_1,\cdots,z_n$ deux \`a
deux distincts, tels que les formes diff\'erentielles $\omega_1:=\frac{\mathrm{d}f_1}{f_1}$ et $\omega_2
:=\frac{\mathrm{d}f_2}{f_2}$ soient de la forme : 

$$\omega_1=\frac{u\mathrm{d}z}{\prod\limits_{i=1}^n(z-x_i)(z-y_i)}
  \;\;\; \mathrm{et}
  \;\;\;\omega_2=\frac{v\mathrm{d}z}{\prod\limits_{i=1}^n(z-x_i)(z-z_i)}$$

Ainsi $\fpp \omega_1+\fpp \omega_2$ est un $L_{m+1,2}$.

R\'eciproquement, tout espace $L_{m+1,n}$ est de cette forme.
\end{theos}

\noindent
\underline{\textit{D\'emonstration}} : 
Vue la forme demand\'ee pour
$\omega_1$, il faut que $f'_1=u$ et donc que $f_1$ soit
de la forme :
$$f_1=(q(z))^2+uz$$
o\`u $q=z^n+q_1z^{n-1}+\cdots +q_n$ est un polyn\^ome de degr\'e
$n$ \`a coefficients dans $k$. De m\^eme, on a $f_2 =(r(z))^2+vz$
o\`u $r =z^n+r_1z^{n-1}+\cdots +r_n$ est un polyn\^ome du m\^eme
type. D\'eterminons donc les polyn\^omes $q$ et $r$.

Remarquons que $f_1(x_i)=(q(x_i))^2+ux_i=0$ ce qui
donne le syst\`eme :
$$\left\{
  \begin{array}{ccc}
  x_1^n+q_1x_1^{n-1}+\cdots+q_n & = & \sqrt{ux_1} \\
  \vdots & & \\ \vdots & & \\
   x_n^n+q_1x_n^{n-1}+\cdots+q_n & = & \sqrt{ux_n}
  \end{array}
  \right.
$$
Vu que les $x_i$ sont distincts, ceci est un syst\`eme de
type Vandermonde, ce qui donne une solution pour les
$q_1,\cdots,q_n$ (et donc pour les $y_1,\cdots,y_n$). De plus, puisque
$f'_1(z)=u$, $f_1$ n'a que des racines simples (donc les
$x_1,\cdots,x_n,y_1,\cdots,y_n$ sont deux \`a deux distincts).
 
On obtient de fa\c con identique que les coefficients du polyn\^ome
$r$ sont obtenus par r\'esolution d'un syst\`eme de Vandermonde. Ceci
fournit les points $z_i$ (et de m\^eme on a que les $x_1,\cdots,x_n,z_1,\cdots,z_n$  sont
deux \`a deux distincts). Il reste \`a v\'erifier que les $y_1,\cdots,y_n,z_1,\cdots,z_n$ 
sont distincts deux \`a deux. 

Soit $\alpha$ une racine commune \`a $f_1$ et $f_2$. Alors :
$$ (q(\alpha ))^2+u\alpha ^{2n-1}= (r(\alpha ))^2+v\alpha ^{2n-1}=0$$
Donc $(vq^2+ur^2)(\alpha )=0=(\sqrt{v}q+\sqrt{u}r)^2(\alpha )$. Or le
polyn\^ome $(\sqrt{v}q+\sqrt{u}r)$ est de degr\'e $n$ et a donc au plus
$n$ racines (qui sont en fait les $x_i$). Finalement les points
$x_1,\cdots,x_n,y_1,\cdots,y_n,z_1,\cdots,z_n$ sont distincts deux \`a deux.

\begin{flushright}
$\square$
\end{flushright}

\subsection{D\'emonstration du r\'esultat principal}

\begin{theos}\label{theo : 2.9}
On consid\`ere le cas $p \geq 3$. 
\begin{enumerate}[1. ]

 \item  Supposons que $m+1=p$. Alors il n'existe pas d'espaces
   vectoriels $L_{m+1,2}$.

 \item Supposons que $m+1=2p$. Alors il existe un espace vectoriel
   $L_{m+1,2}$ si et
   seulement si $p=3$.

 \item Supposons que $m+1=3p$. Alors il n'existe pas d'espaces
   vectoriels $L_{m+1,2}$.
\end{enumerate}
\end{theos}

%\begin{theo}
%\begin{enumerate}
% \item Supposons que $m+1=p$. Alors un tel espace vectoriel n'existe pas.
% \item Supposons que $m+1=2p$. Alors un tel espace vectoriel existe si et
%   seulement si p=3 et dans ce cas, on peut en donner une description.
% \item Supposons que $m+1=3p$. Alors un tel espace vectoriel n'existe pas.
%\end{enumerate}
%\end{theo}

\noindent
\underline{\textit{D\'emonstration}} : Dans les trois cas, la
d\'emonstration se fait par l'absurde et on consid\`erera  donc \`a chaque
fois un espace vectoriel r\'epondant au probl\`eme. Soit $z$ un
param\`etre de $\pp _k-\{\infty\}$ tel que $z=0$ n'est pas p\^ole de
$\omega_1$ et $\omega_2$. On utilise alors dans la d\'emonstration le param\`etre
$x:=\frac{1}{z}$.  

\underline{Le cas $m+1=p$.}

On note toujours $(\omega_1,\omega_2)$ une base d'un espace vectoriel
$L_{m+1,2}$. Ces deux formes s'\'ecrivent :

$$\omega_1=\frac{\mathrm{d}f_1}{f_1}=\sum_{i=0}^p\frac{h'_ix_i}{1-x_ix}\mathrm{d}x,
\;\;\; h'_0=0,\; h'_i\neq 0\; \mathrm{si}\; i \neq 0,\; \sum_i
h'_i\equiv 0 \;\mathrm{mod} \; p$$
$$\omega_2=\frac{\mathrm{d}f_2}{f_2}=\sum_{i=0}^p\frac{h_ix_i}{1-x_ix}\mathrm{d}x,
\;\;\; h_p=0,\; h_i\neq 0 \; \mathrm{si}\; i \neq p,\; \sum_i
h_i\equiv 0 \;\mathrm{mod} \; p$$
et tous les $x_i$ sont distincts.

La forme $\omega_1$ s'\'ecrit aussi :
$$\omega_1=\frac{ux^{p-2}}{\prod_{i=1}^p(1-x_ix)}\mathrm{d}x \;\;\;
\mathrm{avec} \; u\neq 0$$
En identifiant les termes en $x^k$ des d\'eveloppements formels
des deux expressions de $\omega_1$, on trouve que :
$$\sum_{i=1}^ph'_ix_i^k=0 \;\;\; \mathrm{si} \; k\leq p-2\;\;\;\mathrm{et}
\sum_{i=1}^ph'_ix_i^p=u\sum_{i=1}^px_i $$
Or $\sum_{i=1}^ph'_ix_i^p=\sum_{i=1}^p(h'_ix_i)^p $ et vu que $u \neq
0$, il suit que $\sum_{i=1}^px_i=0$. En appliquant 
le m\^eme raisonnement \`a $\omega_2$, il vient que 
$\sum_{i=0}^{p-1}x_i=0$. Ainsi $x_0=x_p$, ce qui fournit la
contradiction attendue. 

%\noindent
%\underline{\textit{Remarque}} :

\vspace{1em}
\underline{Supposons maintenant que $m+1=2p$.}

D'apr\`es le lemme \ref{lemme : 1.5} on a $2p+2$ p\^oles que l'on peut partitionner en $p+1$ couples. On
les note $x_0,y_0,\cdots,x_p,y_p$. Alors, apr\`es renum\'erotation
\'eventuelle, on a que :
\begin{enumerate}[$\bullet$]
 \item $\omega_1+i\omega_2$ a des p\^oles en tous les points sauf en
   $x_i$ et $y_i$ ($i$ varie de $0$ \`a $p-1$).
    
 \item $\omega_p$ a des p\^oles en tous les points sauf en $x_p$ et
   $y_p$.
\end{enumerate}
\noindent
On peut \'ecrire :
$$\omega_1=\sum_{i=0}^p
\left(\frac{h'_ix_i}{1-x_ix}+\frac{k'_iy_i}{1-y_ix}\right)\mathrm{d}x\;\;\;\mathrm{avec}\;h'_0=k'_0=0$$
 $$\omega_2=\sum_{i=0}^p
\left(\frac{h_ix_i}{1-x_ix}+\frac{k_iy_i}{1-y_ix}\right)\mathrm{d}x\;\;\;\mathrm{avec}\;h_p=k_p=0$$

\noindent
\underline{\textit{Etape $1$}} : Montrons que $x_i+y_i$ est une
constante ind\'ependante de $i$.

\vspace{1em}

On pose $s_i=x_i+y_i$ et $p_i=x_iy_i$. Alors $P_i(z)=z^2-s_iz+p_i$;
on a donc, d'apr\`es le paragraphe 1.3, les relations
suivantes :
$$s_i=\frac{as_0+is_p}{a+i}\;\;\;\mathrm{et}\;\;\;p_i=\frac{ap_0+ip_p}{a+i}$$
et $a \notin \fp$.

\noindent
Le m\^eme argument que dans le cas $m+1=p$ montre que :
$$\sum_{i=1}^p(x_i+y_i)=0\;\mathrm{et}\;\sum_{i=0}^{p-1}(x_i+y_i)=0$$
On en d\'eduit donc que $s_0=s_p$, puis finalement que $s_i=$ cste, au
vu de la relation $(a+i)s_i=as_0+is_p$. \finetape

On posera $s_i=s$ dans la suite et $A_i(k)=h_ix_i^k+k_iy_i^k$ pour $k
\geq 0$. 

\vspace{1em}
\noindent
\underline{\textit{Etape $2$}} :Montrons par r\'ecurrence sur $l$ que
:
 $$\sum_{i=0}^{p-1} p_i^lA_i(k)=0\;\;\; 0\leq k \leq 2p-2-2l$$

\vspace{1em}

\begin{enumerate}[$\bullet$]
 \item $l=0$ : La relation annonc\'ee est vraie, car la condition imposant \`a $\omega_2$ d'avoir un z\'ero
d'ordre $(2p-2)$ en z\'ero est :
$$\sum_{i=0}^{p-1} A_i(k)=0, \;\;\; 0\leq k \leq 2p-2$$
 \item Supposons le r\'esultat vrai au rang $l$. On part de l'\'egalit\'e :
$$
\sum_{i=0}^{p-1}p_i^l(A_i(k+2)-sA_i(k+1)+p_iA_i(k))=0\;\;\;\forall
k \geq 0.$$

Alors pour $2\leq k+2\leq 2p-2-2l$ (i.e $0 \leq k\leq 2p-2-2(l+1)$), on a
$$\sum_{i=0}^{p-1}p_i^lA_i(k+2)=0$$ et 
$$\sum_{i=0}^{p-1}sp_i^lA_i(k+1)=0$$
 Donc :
 $$\sum_{i=0}^{p-1} p_i^{l+1}A_i(k)=0\;\;\;\forall k \leq 2p-2-2(l+1)$$ 
\end{enumerate} \nopagebreak \finetape

On aboutit donc \`a :
\begin{eqnarray}
\sum_{i=0}^{p-1} p_i^lA_i(0)  = 0 &\mathrm{pour} & 1\leq l \leq p-1\\
\sum_{i=0}^{p-1} p_i^lA_i(1)  = 0 &\mathrm{pour} & 1\leq l \leq p-2
\end{eqnarray}

\noindent
\underline{\textit{Etape $3$}} : Montrons que $A_i(0)=0$ et qu'il
existe $\beta$ dans $k^*$ tel que $A_i(1)=\beta (a+i)^{p-2}$ pour
$0\leq i \leq p-1$.

\vspace{1em}

On sait que $p_i=p_p+\frac{a(p_0-p_p)}{a+i}=b+\frac{c}{a+i}$ en
posant $b=p_p$ et $c=a(p_0-p_p)\neq 0$. Le syst\`eme (1) implique en
particulier que :
$$\sum_{i=0}^{p-1}\frac{A_i(0)}{(a+i)^l}=0\;\;\mathrm{pour}\;\; 1\leq l \leq
p-1$$
Posons $F(X)=\sum_{i=0}^{p-1}\frac{A_i(0)}{X+i}$. Alors $a$ est
une racine de $F$ d'ordre au moins \'egal \`a $(p-1)$. Puisque
$\sum_{i=0}^{p-1}A_i(0)=0$, le
num\'erateur de $F$ est de degr\'e au plus $(p-2)$, on en
d\'eduit que $F$ est nul et donc que $A_i(0)=0\;$ pour $0\leq i \leq p-1$.

De m\^eme le syst\`eme (2) implique que :
$$\sum_{i=0}^{p-1}\frac{A_i(1)}{(a+i)^l}=0\;\;\mathrm{pour}\;\; 1\leq l \leq
p-2$$ 
Posons $G(X)=\sum_{i=0}^{p-1}\frac{A_i(1)}{(X+i)}$. Alors $a$ est
une racine de $G$ d'ordre au moins \'egale \`a $(p-2)$. Le
num\'erateur de $G$ \'etant de degr\'e au plus $(p-2)$, on a que $G$
est de la forme :
$$G(X)=\sum_{i=0}^{p-1}\frac{A_i(1)}{(X+i)}=\frac{\beta
  (X-a)^{p-2}}{X^p-X} $$o\`u $\beta$ est une constante non nulle (en
effet, $\beta =0$ impliquerait que $A_i(1)=0$ et, puisque $A_i(0)=0$,
on aurait $h_i=k_i=0$). Apr\`es identification des coefficients dans
cette d\'ecomposition en \'el\'ements simples, on aboutit \`a : 
$$A_i(1)=\beta (a+i)^{p-2}$$ \finetape

En r\'esum\'e, on a :
\begin{enumerate}[$\bullet$]
 \item $A_i(0)=0$ donc $k_i=-h_i$.

 \item $A_i(1)=2h_ix_i=\beta (i+a)^{p-2}$

 \item $s_i=0$, apr\`es translation \'eventuelle sur les $x_i,y_i$ (on voit en effet
 que les deux relations pr\'ec\'edentes restent inchang\'ees apr\`es
 translation). En particulier $x_i \neq 0$. 
\end{enumerate}
\noindent
On a donc un syst\`eme :
$$\left\{ \begin{array}{lll}
h_ix_i &=& \frac{\beta}{2}(i+a)^{p-2} \\
-x_i^2 &=& b+\frac{c}{a+i}
\end{array} \right.
x_i \in k, \;0 \leq i \leq p-1$$
Remarquons que ce syst\`eme traduit \`a lui seul les
conditions impos\'ees par le probl\`eme consid\'er\'e.
On calcule :
$$\frac{h_i^2x_i^2}{x_i^2}=-\frac{(\frac{\beta}{2})^2(i+a)^{2(p-2)}}{b+\frac{c}{a+i}}=-\frac{(\frac{\beta}{2})^2(i+a)^{2p-3}}{b(a+i)+c}=h_i^2
$$
donc,
$$1=(-1)^{\frac{p-1}{2}}\frac{(\frac{\beta}{2})^{p-1}(i+a)^{(2p-3)\frac{p-1}{2}}}{(b(a+i)+c)^{\frac{p-1}{2}}}\;\;\;
0\leq i \leq p-1
$$
Ainsi si
$H(X):=(\frac{\beta}{2})^{p-1}(X+a)^{(2p-3)(\frac{p-1}{2})}-(-1)^{\frac{p-1}{2}}(b(X+a)+c)^{\frac{p-1}{2}}$, $H(X)=0$ mod $X^p-X$. En particulier le coefficient de $X^{p-1}$ dans
$H(X)$ modulo $X^p-X$ est nul. Or on a le :

\begin{lemme}\label{lemme : 2.10}
Soit $n$ un entier sup\'erieur ou \'egal \`a $2$ et $p$ un nombre
premier congru \`a $1$ modulo $n$. Alors
le coefficient de $X^{p-1}$ dans $(X+a)^{(np-(n+1))(\frac{p-1}{n})}$ mod
  $X^p-X$ est :$$C_{q}^2
  (a-a^p)^{q-2}$$
avec $q:=\frac{(n-1)p+(n+1)}{n}$.
\end{lemme}

\underline{\textit{D\'emonstration}} :
Remarquons tout d'abord que :
$$(np-(n+1))(\frac{p-1}{n})=p(p-3)+\frac{(n-1)p+n+1}{n}=p(p-3)+q$$
On a donc :
\begin{eqnarray}
(X+a)^{(np-(n+1))(\frac{p-1}{n})} & = & (X^p+a^p)^{p-3}(X+a)^{q} \nonumber
 \\
 & = & (X+a^p)^{p-3}(X+a)^{q} \;\;\;\mathrm{mod}\; (X^p-X) .\nonumber
\end{eqnarray}
 Supposons que $p-1>n$, dans ce cas $q<p$. Notons $T$ le
 coefficient de $X^{p-1}$ dans l' expression $(X+a^p)^{p-3}(X+a)^{q}$, alors :
\begin{eqnarray}
T & = &
\sum_{j=2}^{q}C_{q}^jC_{p-3}^{p-1-j}a^{(q-j)}(a^p)^{(p-3-(p-1-j))}
  \nonumber \\
& = & \sum_{j=2}^{q}C_{q}^jC_{p-3}^{p-1-j}a^{(q-j)}(a^p)^{(j-2)}
  \nonumber \\
& = &
\sum_{j=0}^{(q-2)}C_{q}^{j+2}C_{p-3}^{j}a^{(q-2)-j)}(a^p)^{j} 
\nonumber
\end{eqnarray}
Regardons le terme $C_{p-3}^j$ modulo $p$. On a :
$$C_{p-3}^j \equiv \frac{(-3)(-4)\cdots
  (-(j+2))}{j!} \equiv (-1)^j\frac{(j+1)(j+2)}{2} \equiv (-1)^jC_{j+2}^2$$
Donc :
\begin{eqnarray*}
C_{q}^{j+2}C_{p-3}^{j} & \equiv &
\frac{q(q-1)(q-2)\cdots
  ((q-2)-j+1)}{j!(j+1)(j+2)}\\ & &(-1)^j\frac{(j+1)(j+2)}{2}
 \\
& \equiv & (-1)^j\;\frac{q(q-1)}{2}\\
& &\frac{(q-2)\cdots
  ((q-2)-j+1)}{j!}  \\
& \equiv & (-1)^jC_{q}^2C_{(q-2)}^j 
\end{eqnarray*}
Finalement :
\begin{eqnarray}
T & = &
C_{q}^2\sum_{j=0}^{(q-2)}(-1)^jC_{(q-2)}^{j}a^{((q-2)-j)}(a^p)^{j}
\nonumber \\
& = & C_{q}^2(a-a^p)^{(q-2)} \nonumber
\end{eqnarray}
%et $C_{\frac{p+3}{2}}^2\neq 0$ mod $p$ si $p>3$.

Il reste \`a examiner le cas o\`u $p-1=n$. Dans ce cas $q=p$ et
$$(X+a^p)^{p-3}(X+a)^{q}=(X+a^p)^{p-2} \;\;\;\mathrm{mod}\; (X^p-X).$$
Le coefficient de $X^{p-1}$ dans l' expression
$(X+a^p)^{p-3}(X+a)^{q}$ vaut donc $0$, il co\"incide avec
$C_q^2=C_p^2$ mod $p$.
\begin{flushright}
$\square$
\end{flushright}

Si on regarde ce lemme pour $n=2$, on voit que le coefficient de
$X^{p-1}$ dans $H(X)$ modulo $X^p-X$ est  $\left(\frac{\beta}{2}\right)^{p-1}C_{\frac{p+3}{2}}^{2}
(a^p-a)^{\frac{p-1}{2}}$. Pour $p>3$, on a $C_{\frac{p+3}{2}}^{2} \neq 0 $ et donc
$a^p=a$, ce qui entraine $a \in \fp$ (ce qui est impossible). 

\vspace{1em}
\noindent
\underline{\textit{Remarque}} : Pr\'ecisons ce qui se passe dans le cas $p=3$.

\noindent
Posons $a_1=h_0x_0=\frac{\beta}{2}a$ et
$a_2=h_1x_1-h_0x_0=\frac{\beta}{2}$. Alors
$h_2x_2=\frac{\beta}{2}(a-1)=a_1-a_2$. On sait enfin que
$h'_1x_1+h'_2x_2+h'_3x_3=0$, donc $h'_3x_3=-a_2$ (on utilise le fait
que $h'_1+h_1=0$ et  $h'_2+2h_2=0$). L'ensemble des huit points \{$x_i$,
$y_i$\} est donc l'ensemble :
$$\{\epsilon_1a_1+\epsilon_2a_2,(\epsilon_1,\epsilon_2)\in
\mathbb{F}_3^2 \backslash \{(0,0)\}\}$$

\vspace{1em}
\underline{Supposons enfin que $m+1=3p$.}

On g\'en\'eralise les notations du cas pr\'ec\'edent en prenant
maintenant $x_i,y_i,z_i$ pour les p\^oles et $h_i,k_i,l_i$ les
r\'esidus correspondants. On posera :
\begin{enumerate}[$\bullet$]
 \item $x_i+y_i+z_i=s=$ cste (m\^eme argument que dans le cas
 $m+1=2p$).

 \item $x_iy_i+y_iz_i+x_iz_i=m_i$.

 \item $x_iy_iz_i=p_i$

  \item $A_i(k)=h_ix_i^k+k_iy_i^k+l_iz_i^k$.
\end{enumerate}

\underline{\textit{Etape $1$}} : Montrons que $m_i$ est constant :

\vspace{1em}
On raisonne par l'absurde et on suppose un instant que $m_i$ est non
constant. Cela permet
apr\`es une translation sur les $x_i,y_i,z_i$, de se ramener 
\`a $p_0=p_p$ puis \`a $p_i$ constant (on notera $p_0$ cette
constante).
 
On a encore cette fois-ci les relations :
$$\sum_{i=0}^{p-1} A_i(k)=0 \;\;\; \mathrm{pour} \;0\leq k\leq 3p-2,$$
$$ m_i=\frac{am_0+im_p}{a+i}\;\;\; \mathrm{et}\;\;\;
p_i=\frac{ap_0+ip_p}{a+i}=p_0\neq 0$$

Comme pr\'ec\'edemment, on part de l'\'egalit\'e :
$$A_i(k+3)-sA_i(k+2)+m_iA_i(k+1)-p_0A_i(k)=0 \;\;\;\forall k \geq 0 $$
qui en sommant sur tous les $i$ donne :

$$\sum_{i=0}^{p-1}(A_i(k+3)-sA_i(k+2)+m_iA_i(k+1)-p_0A_i(k))=0 \;\;\;\forall k \geq 0 $$
et donc :
\begin{eqnarray}
\sum_{i=0}^{p-1} m_iA_i(k)=0 \;\;\; \mathrm{pour}\;1\leq k \leq 3p-4
\end{eqnarray}
Il suit de m\^eme pour $k \geq 2$ : 
\begin{eqnarray}
\sum_{i=0}^{p-1}(m_iA_i(k+2)-m_isA_i(k+1)+m_i^2A_i(k)-p_0m_iA_i(k-1))=0
\end{eqnarray}
et donc que :
$$ \sum_{i=0}^{p-1} m_i^2A_i(k)=0 \;\;\;\mathrm{pour} \; 2\leq k \leq
3p-6$$
Une r\'ecurrence comme dans l'\'etape 2 du cas $m+1=2p$ montre que l'on a de plus
g\'en\'eralement :
\begin{eqnarray}
 \sum_{i=0}^{p-1} m_i^lA_i(k)=0 \;\;\;\mathrm{pour} \; l\leq k \leq
3p-2-2l.
\end{eqnarray}
On a alors en particulier que :
$$\sum_{i=0}^{p-1} m_i^lA_i(p-1)=0 \;\;\;\mathrm{si}\;\;\; 1 \leq l
\leq p-1$$
et
$$\sum_{i=0}^{p-1} m_i^lA_i(p)=0 \;\;\;\mathrm{si}\;\;\; 1 \leq l \leq p-1.$$
Sachant que $m_i=\frac{am_0+im_p}{a+i}$, et gr\^ace \`a un argument
analogue \`a celui du cas $m+1=2p$ (i.e on exhibe un polyn\^ome
de degr\'e au plus $p-2$ ayant un z\'ero d'ordre au moins $p-1$), on a
que :
$$ A_i(p-1)=A_i(p)=0$$
Or $A_i(p)=A_i(1)^p$, donc $A_i(1)=0$. L'expression (4)
\'evalu\'ee en $k=1$ donne alors $\sum_{i=0}^{p-1} m_iA_i(0)=0$, i.e, la
relation (3) pour $k=0$. Ceci entra\^ine que la relation (5) est
encore vraie pour $l-1\leq k \leq 3p-2l-2$. On a donc $\sum_{i=0}^{p-1}
m_i^lA_i(p-2)=0 \;\;\mathrm{si}\;\; 1 \leq l \leq p-1$ et donc par la
m\^eme construction $A_i(p-2)=0$.

Finalement, on a $A_i(p-2)=A_i(p-1)=A_i(p)=0$, d'o\`u $h_i=k_i=l_i=0$
par r\'esolution du syst\`eme lin\'eaire, ce qui est absurde. \finetape

On a donc $m_i=$ cste $=m_0$. La m\^eme manipulation que dans le cas
$m+1=2p$ (\'etape 2) donne les relations :
$$\sum_{i=0}^{p-1}A_i(0)p_i^l=0\;\;\;\mathrm{pour}\;1\leq l \leq p-1$$
$$\sum_{i=0}^{p-1}A_i(1)p_i^l=0\;\;\;\mathrm{pour}\;1\leq l \leq p-1$$ 
$$\sum_{i=0}^{p-1}A_i(2)p_i^l=0\;\;\;\mathrm{pour}\;1\leq l \leq p-2$$
et on en tire de la m\^eme fa\c con que $A_i(0)=A_i(1)=0$ et $A_i(2)$
est de la forme $\beta (i+a)^{p-2}$(l'argument est le m\^eme : on
\'ecrit qu'un polyn\^ome de degr\'e au plus $p-2$ a une racine d'ordre
$p-1$ ou $p-2$ selon les cas). On distingue alors deux cas :

\vspace{1em}
\noindent
\underline{$1$er cas } : $p=3$

Puisque $A_i(0)=h_i+k_i+l_i=0$ on a $h_i=k_i=l_i=\pm
1=\epsilon_i$). On obtient $A_i(2)= \beta (i+a)=\epsilon_i(x_i^2+y_i^2+z_i^2)$, et
$A_i(1)=0=x_i+y_i+z_i$. D'o\`u \nolinebreak:
\begin{eqnarray}
(x_i+y_i+z_i)^2 & = & x_i^2+y_i^2+z_i^2 +2m_0 \nonumber \\
0&=&\epsilon_i \beta(i+a) +2m_0 \nonumber
\end{eqnarray} 
c'est-\`a-dire $\epsilon_i(i+a)$ est une constante. Il existe au moins
deux valeurs de $\epsilon_i$ \'egales ce qui donne la contradiction attendue.

\vspace{1em}
\noindent
\underline{$2$\`eme cas} : $p \neq 3$. On se ram\`ene \`a $s=0$ (par
translation). On a :
\begin{eqnarray}
\left(\begin{array}{ccc}1 & 1 &1 \\ 
                        x_i & y_i & z_i \\
                        x_i^2 & y_i^2 & z_i^2
\end{array} \right)
 \left(\begin{array}{ccc} h_i \\ k_i \\ l_i  \end{array}\right)
& = &  \left(\begin{array}{ccc} 0 \\ 0 \\ \beta (i+a)^{p-2}
  \end{array}\right) \nonumber 
\end{eqnarray}
et donc,
\begin{eqnarray}
h_i & = & \frac{\beta(i+a)^{p-2}}{\Delta}(z_i-y_i) \nonumber \\
k_i & = & \frac{\beta(i+a)^{p-2}}{\Delta}(x_i-z_i) \nonumber \\
l_i & = & \frac{\beta(i+a)^{p-2}}{\Delta}(y_i-x_i) \nonumber 
\end{eqnarray}
\noindent
o\`u $\Delta$ d\'esigne le d\'eterminant de Vandermonde de la matrice
\'ecrite plus haut. Il suit que :
$$h_ik_il_i=\frac{\beta ^3 (i+a)^{3(p-2)}}{\Delta ^2}$$
et 
\begin{eqnarray}
h_ik_i+h_il_i+k_il_i & = & \frac{\beta ^2 (i+a)^{2(p-2)}}{\Delta
    ^2}((z_i-y_i)(x_i-z_i)+(z_i-y_i)(y_i-x_i)\nonumber \\ & + &(x_i-z_i)(y_i-x_i))
    \nonumber \\
  & = & \frac{\beta ^2 (i+a)^{2(p-2)}}{\Delta ^2}(m_0-(x_i^2+y_i^2+z_i^2)) \nonumber \\
  & = & 3\frac{\beta ^2 (i+a)^{2(p-2)}}{\Delta ^2}m_0 \nonumber
\end{eqnarray}

\noindent
\underline{$1$er sous-cas} : $m_0\neq 0$. Alors
$\frac{h_ik_il_i}{h_ik_i+h_il_i+k_il_i}=\frac{1}{3m_0}\beta(i+a)^{p-2}\;\in
\fp$. Donc $(\frac{i}{a}+1)^{p-2} \in \fp$. Posons
$A=\frac{1}{a}$, alors $(Ai+1)^{p(p-2)}-(Ai+1)^{p-2}=0$ $\forall i \in
\fp$. Posons $F(X)=(AX+1)^{p(p-2)}-(AX+1)^{p-2}$, alors :
$$F(X)=(A^pX+1)^{p-2}-(AX+1)^{p-2}=0\;\mathrm{mod}\;(X^p-X)$$
D'o\`u $A^p=A$ et donc $a\in \fp$ ce qui est absurde.

\vspace{1em}
\noindent
\underline{$2$\`eme sous-cas} : $m_0=0$. Dans ce
cas $x_i^3=p_i=b+\frac{c}{a+i}$, avec $b=p_p$ et $c=a(p_0-p_p)$. On a
la m\^eme relation pour $y_i$ et $z_i$, donc $y_i=jx_i$, $z_i=j^2x_i$,
avec $j^3=1$, $j \neq 1$ (quitte \`a \'echanger $y_i$ et $z_i$ on peut
supposer que $j$ ne d\'epend pas de $i$).

$\bullet$ Si $j\notin \fp$, alors :
\begin{eqnarray}
 h_ix_i+k_iy_i+l_iz_i & = & 0 \nonumber \\
 x_i(h_i+k_ij-(1+j)l_i) & = & 0 \nonumber \\
 (h_i-l_i)+j(k_i-l_i) & = & 0 \nonumber
\end{eqnarray} 
Donc $h_i=k_i=l_i$. Comme $A_i(0)=h_i+k_i+l_i=0$ et que $p\neq 3$, on
obtient une absurdit\'e.

$\bullet$ Donc $j \in \fp$ et $p \equiv 1$ mod $3$. Des \'egalit\'es $A_i(0)=A_i(1)=0$
 , on tire le syst\`eme lin\'eaire en $h_i,k_i,l_i$ :
$$\left\{\begin{array}{ccc} h_i+k_i+l_i & = & 0 \\
                           h_i+jk_i-(1+j)l_i & = & 0
    \end{array}\right.
$$
ce qui permet par exemple d'exprimer $k_i$ et $l_i$ en fonction de
$h_i$ :
$$\left\{\begin{array}{ccccc} k_i & = & \frac{-(2+j)}{2j+1} h_i & = & \mu
      h_i \\
 l_i & = & \frac{1-j}{2j+1}h_i & = & \lambda h_i
\end{array}\right.
$$
$\mu$, $\lambda \in \fp$, (ind\'ependants de $i$). Donc :
\begin{eqnarray}
  h_ix_i^2+k_iy_i^2+l_iz_i^2 & = & \beta (i+a)^{p-2}  =  A_i(2)
 \nonumber \\
  h_ix_i^2(1+\mu j^2+\lambda j^4)  & = & \beta (i+a)^{p-2} 
 \nonumber \\
  h_ix_i^2 & = &  \beta ' (i+a)^{p-2}  \nonumber
\end{eqnarray}
en posant $\beta '=\beta(1+\mu j^2+\lambda j^4)^{-1}$. Puisque
$x_i^3=b+\frac{c}{a+i}$, il suit que :

$$\frac{h_i^3x_i^6}{x_i^6}=\frac{{\beta
  '}^3(i+a)^{3(p-2)}}{(b+\frac{c}{a+i})^2}=\frac{{\beta '}^3(i+a)^{3p-4}}{(b(a+i)+c)^2}=h_i^3$$
donc,
$$1=\frac{{\beta '}^{p-1}(i+a)^{(3p-4)\frac{p-1}{3}}}{(b(a+i)+c)^{2\frac{p-1}{3}}}$$

Posons
$G(X)={\beta '}^{p-1}(X+a)^{(3p-4)(\frac{p-1}{3})}-(b(X+a)+c)^{2\frac{p-1}{3}}$.
Alors $G(X)=0$ mod $X^p-X$. En particulier le coefficient de $X^{p-1}$ dans
$G(X)$ modulo $X^p-X$ est nul. 

On peut appliquer le lemme \ref{lemme : 2.10} pour $n=3$, (notons que $p \equiv 1$
mod $3$); le coefficient
en $X^{p-1}$ de $G(X)$ modulo $X^p-X$ est $\beta^{'p-1}C_{\frac{2p+4}{3}}^{2}
(a^p-a)^{\frac{2(p-1)}{3}}$. Or, $C_{\frac{2p+4}{3}}^{2} \neq 0 $ donc $a^p=a$,
et donc $a \in \fp$, d'o\`u la contradiction. 
\finetape

Dans tous les cas, il n'y a pas de $L_{3p,2}$ pour $p>2$.
% \begin{flushright}
%$\square$
%\end{flushright}

 \subsection{Exemples d'espaces vectoriels $L_{m+1,n}$}

%Dans \cite{mat}, on a des exemples d'espaces $L_{m+1,n}$ qui
%permettent de construire des actions de $(\Z /p\Z)^n$ sur le disque
%ouvert $p$-adique. Explicitons ce que l'on trouve \`a
%partir de \cite{mat}. 
On peut expliquer la construction qui se trouve dans \cite{mat}
d'actions de $(\Z/p\Z)^n$ sur le disque ouvert $p$-adique par la
pr\'esence cach\'ee d'espaces $L_{m+1,n}$. Explicitons cela.

%L'id\'ee directrice de cette construction part du constat que la forme
%diff\'erentielle $\omega:=\frac{u\mathrm{d}z}{z^{p-1}-\alpha}$ est
%logarithmique ($\alpha \in k^*$ et $u=\frac{\alpha}{x_i}$, o\`u $x_i$
%est une des racines du polyn\^ome $z^{p-1}-\alpha$). On constate alors
%que si l'on fait le changement de variables $z=P(t)$ o\`u $P$ est un
%polyn\^ome de d\'eriv\'e constant, la nouvelle forme obtenue est
%encore logarithmique et a un seul z\'ero \`a l'infini. 

Dans cette construction, on utilise le fait que la forme
$\omega:=\frac{u\mathrm{d}z}{z^{p-1}-\alpha}$ ($\alpha \in k^*$ et
$u=\frac{\alpha}{x_i}$, o\`u $x_i$ est une des racines du polyn\^ome
$z^{p-1}-\alpha$) est logarithmique, et la remarque 4 du paragraphe 1.3.

On se donne un entier $n$ sup\'erieur ou \'egal \`a $2$. Consid\'erons
les formes diff\'erentielles suivantes :
$$\omega_j=\frac{u_j\mathrm{d}z}{\prod\limits_{\substack{(\epsilon_1,\cdots,\epsilon_n)\in
   \{0,\cdots,p-1\}^n \\ \epsilon_j
   \neq0}}(z-\sum\limits_{i=1}^n\epsilon_ia_i)}$$
o\`u $u_j$ est une constante que l'on va montrer pouvoir choisir ``convenablement''
pour que la forme $\omega_j$ soit logarithmique.
%De m\^eme :
%$$\omega_2=\frac{v\mathrm{d}x}{\prod\limits_{\substack{(\epsilon_1,\cdots,\epsilon_n)\in
%   \{0,\cdots,p-1\}^n \\ \epsilon_2
%   \neq0}}(x+\sum\limits_{i=1}^n\epsilon_ia_i)}$$
%($u$ et $v$  sont des constantes que l'on va expliciter).
%V\'erifions que ces formes diff\'erentielles conviennent. Le fait que
%$\omega_1$ et $\omega_2$ engendrent un espace vectoriel de formes
%diff\'erentielles avec les bonnes conditions de z\'eros et de p\^oles
%est assez \'evident.
On a :
\begin{eqnarray}
  \omega_1 & = & \frac{u_1\mathrm{d}z}{\prod\limits_{\substack{(\epsilon_1,\cdots,\epsilon_n)\in
   \{0,\cdots,p-1\}^n \\ \epsilon_1
   \neq0}}(z-\sum\limits_{i=1}^n\epsilon_ia_i)} \nonumber \\
& = & \frac{u_1\mathrm{d}z}{\prod\limits_{j=1}^{p-1}\prod\limits_{(\epsilon_2,\cdots,\epsilon_n)\in
   \{0,\cdots,p-1\}^n}(z-ja_1+\sum\limits_{i=2}^n\epsilon_ia_i)}
   \nonumber 
\end{eqnarray}
Notons 
$$Ad_1(z)=\prod\limits_{(\epsilon_2,\cdots,\epsilon_n)\in
   \{0,\cdots,p-1\}^n}(z-\sum\limits_{i=2}^n\epsilon_ia_i)$$
Alors $Ad_1$ est un polyn\^ome additif; $\omega_1$ s'\'ecrit alors :
\begin{eqnarray}
  \omega_1 & = &
  \frac{u_1\mathrm{d}z}{\prod\limits_{j=1}^{p-1}Ad_1(z-ja_1)} \nonumber \\
 & = &
  \frac{u_1\mathrm{d}z}{\prod\limits_{j=1}^{p-1}(Ad_1(z)-jAd_1(a_1))} \nonumber \\
 & = & \frac{u_1\mathrm{d}z}{Ad_1(z)^{p-1}-Ad_1(a_1)^{p-1}} \nonumber
\end{eqnarray}
On peut \'ecrire $Ad_1$ sous la forme $\alpha_1 z +P_1(z^p)$, car
$Ad_1$ est additif (cf remarque 4 du paragraphe 1.3). En
particulier, $Ad_1\,'(z)=\alpha_1$. Posons $Q(z)=Ad_1(z)^{p-1}-Ad_1(a_1)^{p-1}$
et calculons $Q'(\sum_{i=1}^n\epsilon_ia_i)$ pour $\epsilon_i \in
\{0,\cdots,p-1\}$, $\epsilon_1 \neq 0$.
\begin{eqnarray}
 Q'(\sum_{i=1}^n\epsilon_ia_i) & = & -\alpha_1
 Ad_1(\sum_{i=1}^n\epsilon_ia_i)^{p-2} \nonumber \\
& = & -\alpha_1 (\sum_{i=1}^n \epsilon_i Ad_1(a_i))^{p-2} \nonumber \\
& = & -\alpha_1 (\epsilon_1 Ad_1(a_1))^{p-2} \nonumber
\end{eqnarray}
Posons alors $u_1=-\alpha_1 Ad_1(a_1)^{p-2}$. Alors :
\begin{eqnarray}
\omega_1= \frac{u_1\mathrm{d}z}{Q(z)} & = & \sum_{\substack{(\epsilon_1,\cdots,\epsilon_n)\in
   \{0,\cdots,p-1\}^n \\ \epsilon_1 \neq 0}}
   \frac{\frac{u_1\mathrm{d}z}{Q'(\sum_{i=1}^n\epsilon_ia_i)}}{(z-\sum\limits_{i=1}^n\epsilon_ia_i)}\nonumber
   \\ & = &  \sum_{\substack{(\epsilon_1,\cdots,\epsilon_n)\in
   \{0,\cdots,p-1\}^n \\ \epsilon_1 \neq 0}} \frac{\epsilon_1\mathrm{d}z}{(z-\sum\limits_{i=1}^n\epsilon_ia_i)}\nonumber
\end{eqnarray}
ce qui prouve que $\omega_1$ est bien logarithmique. De m\^eme, on
peut trouver $u_j$ pour que $\omega_j$ soit logarithmique; $\omega_j$
s'\'ecrit alors :
$$\omega_j=  \sum_{\substack{(\epsilon_1,\cdots,\epsilon_n)\in
   \{0,\cdots,p-1\}^n \\ \epsilon_j \neq 0}}
\frac{\epsilon_j\mathrm{d}z}{(z-\sum\limits_{i=1}^n\epsilon_ia_i)}$$
%Une fois que l'on a les $\omega_j$ sous cette forme,
%il devient alors clair que $\forall [a_1,\cdots,a_n] \in
%\mathbb{P}_{\fp}^{n-1},\; a_1\omega_1+\cdots +a_n\omega_n$ a le bon nombre de
%p\^oles (i.e $ a_1\omega_1+\cdots +a_n\omega_n$ a des p\^oles en tous les
%$-\sum_{i=1}^n\epsilon_ia_i$ tels que $a_1\epsilon_1+\cdots +a_n\epsilon_n \neq
%0$). Le $\fp$-espace vectoriel engendr\'e par les $\omega_j$ est alors
%un $L_{m+1,n}$.

Des consid\'erations de degr\'e montrent que le d\'eterminant de Moore

\noindent
$\Delta(u_1,\cdots,u_n)$ est un polyn\^ome en les $(a_i)_{1\leq
  i\leq n}$ non nul. Ainsi si $(a_1,\cdots,a_n)$

\noindent
$ \in
k^n-V(\Delta(u_1,\cdots,u_n))$, alors $(u_1 \cdots ,u_n)$ sont
$\fp$-lin\'eairement ind\'ependants (c'est la condition (*) de
\cite{mat}).

Sous cette derni\`ere condition, montrons que $<\omega_1,\cdots,\omega_n>$ est un
$L_{m+1,n}$. Puisque $\Delta(a_1,\cdots,a_n)=\Delta(a_1,\cdots,
a_{n-1})Ad_n(a_n)$ et que
$u_n=-\Delta(a_1,\cdots,a_{n-1})^{p-1}Ad_n(a_n)^{p-2}\neq 0$, il suit
que $a_1,\cdots, a_n$ sont $\fp$-lin\'eaire\-ment ind\'ependants.

Soit $(b_1,\cdots,b_n) \in \fp^n-\{0\}$; alors
$b_1\omega_1+\cdots + b_n\omega_n$ a un z\'ero d'ordre $m-1$ \`a
l'infini et $m+1=p^n-p^{n-1}$ p\^oles qui sont les :
$$\{\sum_{i=1}^n \epsilon_ia_i,\; \mathrm{avec}
\;b_1\epsilon_1+\cdots +b_n \epsilon_n \neq 0\}.$$

En r\'esum\'e, si $(a_1,\cdots,a_n)$ v\'erifie la condition (*) de
\cite{mat}, on d\'efinit :
$$\omega_j:=\frac{u_j\mathrm{d}z}{\prod\limits_{\substack{(\epsilon_1,\cdots,\epsilon_n)\in
   \{0,\cdots,p-1\}^n \\ \epsilon_j
   \neq0}}(z-\sum\limits_{i=1}^n\epsilon_ia_i)}$$
Alors $<\omega_1,\cdots,\omega_n>$ est un $L_{m+1,n}$.

\vspace{1em}
\noindent
\underline{\textit{Remarque 1}} : Si on reprend les arguments de la
remarque 4 du paragraphe 1.3, on peut construire par changement
de variables d'autres exemples d'espaces $L_{m+1,n}$.

\vspace{1em}
\noindent
\underline{\textit{Remarque 2}} : Pour chaque exemple d'espaces
$L_{m+1,n}$ ainsi construits, on constate que $m+1$ est un multiple de
$p^{n-1}(p-1)$ . Il est tentant de penser (cf. th\'eor\`eme \ref{theo :  2.9})
que cette condition est n\'ecessaire.

%Soit $\mathrm{G}=(\Z/p\Z)^2$, et $\mathrm{R}$ un anneau de valuation
%discr\`ete dominant l'anneau des vecteurs de Witt de $k$. On se donne
%une extension $k[[z]]/k[[z]]^{\mathrm{G}}$, et on veut la relever en
%$\mathrm{R}[[Z]]/\mathrm{R}[[Z]]^{\mathrm{G}}$. On
%regarde le cas particulier o\`u chacune des sous-extensions de $
%k[[z]]^{\mathrm{G}}$ d'ordre $p$ a un conducteur \'egal \`a $m+1=p$. 
\subsection{Applications}

Le probl\`eme de l'existence de ces espaces $L_{m+1,n}$ est intimement
li\'e aux actions de $(\Z/p\Z)^n$ sur le disque ouvert $p$-adique.

\subsubsection{Action de $(\Z/p\Z)^n$ sur le disque ouvert $p$-adique}

(Pour plus de pr\'ecisions sur les rappels qui vont suivre, on renvoie
\`a \cite{gre1} et \cite{gre2}.) 

Soit $\mathrm{R}$ un anneau de valuation
discr\`ete dominant l'anneau des vecteurs de Witt de $k$. 
Soit $\mathrm{D}_0=\mathrm{Spec}(\mathrm{R}[[Z]])$ et
$\sigma$ un automorphisme d'ordre $p$ agissant sur $\mathrm{D}_0$ et
ayant $m+1$ points fixes. On note $\dd$ le mod\`ele semi-stable
minimal qui d\'eploie les $m+1$ points fixes en des points lisses et distincts dans la fibre sp\'eciale
$\ddd$. La fibre sp\'eciale est alors un arbre de droites projectives,
on montre que les sp\'ecialisations des points  fixes se trouvent
dans les composantes terminales de l'arbre. Notons $\dd
':=\dd/<\sigma>$. Les fibres sp\'eciales $\ddd$ et $\ddd '$ sont alors
hom\'eomorphes via le morphisme de passage au quotient par $\sigma$.

Consid\'erons une composante terminale $E'$ de $\ddd'$, alors
des espaces $L_{m+1,1}$ apparaissent quand on analyse
la d\'eg\'enerescence du $\mu _p$-torseur induit par $\sigma$
sur le disque ferm\'e correspondant \`a la composante
$E'$. Pr\'ecis\'ement, si on note $x_0, \cdots, x_m$ les points
fixes de $\sigma$ qui se sp\'ecialisent dans la composante $E'$, on montre qu'il
existe $f \in k(E')$ telle que $\mathrm{ord}_{\infty} f=0$ et
$\mathrm{d}f$ a son diviseur \`a support dans $\{x_i\}_i \cup
\{\infty\}$, $\mathrm{ord}_{x_i}\mathrm{d}f=h_i-1\;(\mathrm{mod}\;
p),\mathrm{ord}_{\infty}\mathrm{d}f=m-1$ (et $\sum h_i=0$),(cf
\cite{gre2} th\'eor\`eme III.3.1).  
%Vu que $m<p$, $\ddd$ est une droite projective sur
%laquelle on a $m+1$ points $x_0,\cdots,x_m$ \'equidistants qui
%correspondent aux sp\'ecialisations des points fixes par $\sigma$ (cf
%\cite{gre2} th\'eor\`eme III.3.1).
Apr\`es un changement de param\`etre, on voit que $f$ est de la forme :
$$\prod_{i=0}^m(1-x_ix)^{h_i}$$ 
$$ \mathrm{et}\;
\frac{\mathrm{d}f}{f}=-\sum_{i=0}^m\frac{h_ix_i}{1-x_ix}\mathrm{d}x=\frac{-ux^{m-1}}{\prod_{i=0}^m(1-x_ix)},\;u\in
k^*.$$  

La g\'eom\'etrie la plus simple qui peut intervenir est la
g\'eom\'etrie \'equidistante, i.e la fibre sp\'eciale $\ddd$ est
r\'eduite \`a une droite projective. Dans la cas o\`u $\sigma \neq Id$
mod $\mathfrak{m}$ (o\`u $\mathfrak{m}$ est l'id\'eal maximal de R), la
distance mutuelle entre les points fixes est
$|\zeta-1|^{\frac{1}{m}}$; elle est donc d\'etermin\'ee par le
conducteur de l'extension. 

Plus g\'en\'eralement, on peut d\'efinir des donn\'ees combinatoires
et diff\'erentielles sur les autres composantes
(cf. \cite{gre2}). Henrio a \'etabli dans \cite{hen} la r\'eciproque,
c'est-\`a-dire, reconstruire un automorphisme d'ordre $p$ \`a partir
de ces donn\'ees.

Nous nous proposons dans ce qui suit d'aborder sous le m\^eme angle
l'action du groupe $G:=(\Z/p\Z)^n$. Dans le cas d'une action de $G$ sur le disque ouvert
$p$-adique, la description des donn\'ees combinatoires et
diff\'erentielles est plus d\'elicate. Nous allons examiner ici le cas
de la combinatoire la plus simple, i.e. le cas o\`u le lieu de
branchement du $G$-torseur correspondant est \'equidistant.

On se donne donc un $G$-torseur au dessus de SpecR$[[T]]:=$
SpecR$[[Z]]^G$. On a ainsi $n$ rev\^etements $p$-cycliques
SpecR$[[Z_i]] \rightarrow$ SpecR$[[T]]$ donn\'es par les \'equations
$Z_i^p=f_i(T)$ (o\`u $f_i \in $R$[T]$). Consid\'erons une extension de
SpecR$[[T]]$ $p$-cyclique interm\'ediaire; elle est donn\'ee par une
\'equation du type $Y^p=f_1^{\epsilon_1} \cdots f_n^{\epsilon_n}$ avec
$(\epsilon_1, \cdots ,\epsilon_n) \in \fp^n-\{0\}$. L'hypoth\`ese
faite sur le lieu de branchement (g\'eom\'etrie \'equidistante) impose
alors que ce rev\^etement a pour conducteur $m+1$. En particulier la
forme diff\'erentielle logarithmique  associ\'ee \`a ce rev\^etement a
$m+1$ p\^oles distincts et un z\'ero d'ordre $m-1$ \`a l'infini. On
obtient donc ainsi un espace $L_{m+1,n}$. 

Nous remarquons aussi qu'une action de $(\Z/p\Z)^n$ sur R$[[Z]]$ induit
en r\'eduction (i.e modulo $\mathfrak{m}$) une action de
$(\Z/p\Z)^n$ sur $k[[z]]$ (qui peut \^etre triviale).

%$W(k)$
%\noindent
%\underline{Question 1} : Inversement, on peut se demander si \`a partir d'un espace
%$L_{m+1,n}$, on peut fabriquer une action de $(\Z/p\Z)^n$ sur le
%disque ouvert $p$-adique. En d'autres termes, \'etant donn\'e un
%$\mu_p^n$-torseur au dessus de $\pp_k$, peut-on le d\'eformer en une
%action de $(\Z/p\Z)^n$ sur le disque ouvert $p$-adique. Le paragraphe
%suivant apporte une r\'eponse positive \`a cette question. \finetape

%Partons encore une fois d'une action de $(\Z/p\Z)^n$ sur R$[[T]]$. Si
%on regarde ceci en r\'eduction (i.e modulo $\mathfrak{m}$, o\`u
%$\mathfrak{m}$ est l'id\'eal maximal de R), on trouve une action de
%$(\Z/p\Z)^n$ sur $k[[t]]$. 

%\noindent
%\underline{Question 2} : Etant donn\`e un $(\Z/p\Z)^n$-torseur au dessus de Spec$k[[t]]$,
%peut-on le d\'eformer en un $(\Z/p\Z)^n$-torseur au dessus de
%Spec$R[[T]]$ o\`u $R$ est un anneau de valuation discr\`ete dominant
%$W(k)$ ? 

%On apportera une r\'eponse \`a ce probl\`eme dans deux cas particuliers.

%\noindent
%\underline{\textit{ Remarque}} : Dans \cite{hen}, Henrio donne un r\'esultat
%bien plus g\'en\'eral : il montre que l'on peut construire des
%R-automorphismes d'ordre $p$ de R$[[Z]]$ si l'on dispose au
%pr\'ealable d'un arbre d'Hurwitz assujeti \`a certaines conditions
%suppl\'ementaires. (Notre r\'esultat correspond \`a l'arbre le
%plus simple).
\subsubsection{Construction de $(\Z/p\Z)^n$-torseurs \`a partir
  d'espaces $L_{m+1,n}$ }

Dans un premier temps , nous allons montrer qu'un espace $L_{m+1,n}$
donne naissance \`a une action de $(\Z/p\Z)^n$ sur le disque ouvert
$p$-adique (nous pr\'ecisons \'egalement le $(\Z/p\Z)^n$-torseur obtenu
en r\'eduction modulo $\mathfrak{m}$).

%Ensuite, nous examinons la question du rel\`evement d'un
%$(\Z/p\Z)^n$-torseur au dessus de Spec$k[[t]]$  en un
%$(\Z/p\Z)^n$-torseur au dessus de Spec$R[[T]]$. Nous apporterons une
%r\'eponse \`a ce probl\`eme dans deux cas particuliers. 

% Nous allons montrer que si l'on se donne un espace $L_{m+1,n}$,
% on peut en d\'eduire une r\'ealisation du groupe $G:=(\Z /p\Z)^n$ comme
% groupe d'automorphis\-mes du disque ouvert $p$-adique.

Plus pr\'ecis\'ement, on a le th\'eor\`eme suivant :
\begin{theos}\label{theo : 2.11}
On consid\`ere un $L_{m+1,n}$ et une base $\omega_1,\cdots,\omega_n$
de cet espace, chaque $\omega_i$ s'\'ecrivant
$\frac{\mathrm{d}f_i}{f_i}$. Soit $\zeta$ une racine primitive $p$-i\`eme de
l'unit\'e et R=W$(k)[\pi]$ o\`u $\pi ^m:=\lambda:=\zeta -1$, on note K=Frac(R). Alors on peut
trouver $F_i \in $R$[X]$ relevant $f_i$ tels que le produit fibr\'e
des rev\^etements de $\pp_{\mathrm{K}}$ donn\'es par les \'equations
$Y_i^p=F_i(X)$ induisent apr\`es normalisation un rev\^etement de
$\pp_{\mathrm{K}}$ galoisien de groupe $(\Z/p\Z)^n$ ayant bonne
r\'eduction relativement \`a la valuation de Gauss $T:=\pi^{-p}X$. La
fibre sp\'eciale du mod\`ele lisse correspondant est un rev\^etement  \'etale, galoisien de
groupe  $(\Z/p\Z)^n$ de la droite affine $\mathbb{A}_k^1$.
\end{theos} 

La d\'emonstration suit les m\'ethodes utilis\'ees
dans \cite{mat}. Nous allons l'adapter au cas qui nous pr\'eoccupe. 

Nous montrons d'abord le lemme suivant :
\begin{lemme}\label{lemme : 2.12}
%On conserve les m\^emes notations. Soit $g=\prod_{i=0}^m
%(1-x_ix)^{h_i}$ telle que $\frac{\mathrm{d}g}{g}$ soit de la forme
%$\frac{ux^{m-1}\mathrm{d}x}{\prod_{i=0}^m (1-x_ix)}$ ( en d'autres
%termes $\frac{\mathrm{d}g}{g}$ engendre un $L_{m+1,1}$). Soit $X_i \in
%W(k)$ des rel\`evements de $x_i$; on pose $G(X):=\prod_{i=0}^m
%(1-X_iX)^{h_i}$. Alors il existe
%$\hat{Q}(X),\hat{R}(X),\hat{S}(X) \in W(k)[X]$ et $U\in W(k)^*$ tels que :
% $$G=(1+X\hat{Q}(X))^p+UX^m
% ( 1+X\hat{R}(X))+p\hat{S}(X)$$
Soit $\omega_1,\cdots,\omega_n$ une base d'un espace $L_{m+1,n}$; soit
$(x_i)_{1\leq i\leq T}$ la r\'eunion des p\^oles de $\omega_j$ pour
$1\leq j \leq n$ et $(x_i)_{i\in I_j}$ les p\^oles de
$\omega_j$. Chaque $\omega_j$ s'\'ecrit $\frac{\mathrm{d}f_j}{f_j}$
avec $f_j=\prod_{i=1}^T(1-x_ix)^{h_{ij}}$ et $h_{ij}=0$ pour $i
\notin I_j$. Soit $X_i \in W(k)$ des rel\`evements de $x_i$ pour
$1\leq i \leq T$ . On pose $F_j(X):= \prod_{i=1}^T
(1-X_iX)^{h_{ij}}$. Alors il existe
$\hat{Q}_j(X),\hat{R}_j(X),\hat{S}_j(X) \in W(k)[X]$ et $U_j\in W(k)$ inversible tels que :
 $$F_j(X)=(1+X\hat{Q}_j(X))^p+U_jX^m
 ( 1+X\hat{R}_j(X))+p\hat{S}_j(X)  \hspace{3em} (*)$$

\end{lemme}

\noindent
\underline{\textit{D\'emonstration}} : On a :
$$f'_j=\frac{u_jx^{m-1}}{\prod_{i=1}^T (1-x_ix)}\prod_{i=1}^T
(1-x_ix)^{h_{ij}}=u_jx^{m-1}(1+xr(x))$$
o\`u $r(x)$ est un polyn\^ome dans lequel on a regroup\'e tous les
termes de degr\'e sup\'erieur. Le polyn\^ome $f_j$ est donc de la forme
:
$$f_j=(1+xq(x))^p+\frac{u_jx^m}{m}(1+x\tilde{r}(x)).$$ Donc $F_j$ qui est un
rel\`evement de $g$ s'\'ecrit :
 $$F_j=(1+X\hat{Q}_j(X))^p+U_jX^m( 1+X\hat{R}_j(X))+p\hat{S}_j(X)$$
\begin{flushright}
$\square$
\end{flushright}

\noindent
\underline{\textit{D\'emonstration du th\'eor\`eme}} :
L'approximation (*) du lemme \ref{lemme : 2.12} n'est a priori pas suffisante pour
garantir que les $F_j$ satisfassent le th\'eor\`eme. On va am\'eliorer
cette approximation en utilisant l'automorphisme de Frobenius. L'action
du Frobenius inverse sur $k[t]$ est d\'efinie de la fa\c
con suivante \nolinebreak: si $f:=\sum a_ix^i \in k[t]$ alors on pose $f^{F^{-1}}:=\sum
a_i^{\frac{1}{p}}x^i$.  Cette op\'eration commute avec la d\'erivation (i.e
$(f^{F^{-1}})'=(f')^{F^{-1}}$). On peut donc \'etendre cette action aux formes
diff\'erentielles que l'on consid\`ere. En particulier, si on a un
espace $L_{m+1,n}$ engendr\'e par les $n$ formes diff\'erentielles $
\omega_1, \cdots, \omega_n$ alors on en d\'eduit que le $\fp$-espace
vectoriel engendr\'e par les formes $\omega_1^{F^{-1}}, \cdots, \omega_n^{F^{-1}}$
est encore un espace $L_{m+1,n}$. 

%On consid\`ere un espace $L_{m+1,n}$ de base
%$\frac{\mathrm{d}f_1}{f_1},\cdots,\frac{\mathrm{d}f_n}{f_n}$.
On choisit une des fonctions $f_j$ (que l'on appelle $f$ dans la suite
pour ne pas surcharger les notations; de la m\^eme fa\c con on notera
$h_i$ \`a la place de $h_{ij}$). Nous allons montrer qu'il existe des $X_i$ relevant
$x_i$ tels que la fonction $F$ d\'efinie par $F:=\prod_{i=1}^T
(1-X_iX)^{h_i}$  soit de la forme :

$$F(X)=(1+XQ(X))^p+U^pX^m(1+XR(X))+pX^{\frac{(m+1)}{p}}S(X)+p^2T(X)$$
avec $Q(X),R(X),S(X),T(X) \in W(k)[X]$ et $U \in W(k)$ inversible.

%Nous allons montrer le lemme suivant :
%\begin{lemme}
%Soit $\omega$ une forme diff\'erentielle non nulle appartenant \`a un
%espace $L_{m+1,n}$. $\omega$ s'\'ecrit $\frac{\mathrm{d}f}{f}$ avec
%$f$ de la forme :
%$$f:=\prod_{i=0}^m(1-x_ix)^{h_i}$$
%apr\`es un choix convenable de param\`etre.

%Alors il existe $F \in W(k)[X]$ relevant $f$ tel que :
%$$F(X)=(1+XQ(X))^p+U^pX^m(1+XR(X))+pX^{\frac{(m+1)}{p}}S(X)+p^2T(X)$$
%avec $Q(X),R(X),S(X),T(X) \in W(k)[X]$ et $U \in W(k)^*$  
%\end{lemme}

%\noindent
%\underline{\textit{D\'emonstration}} : L'id\'ee de la d\'emonstration
%est d'utiliser la structure de Frobenius.

Soit $y_i \in k$ tels que
$y_i^p=x_i$; prenons $Y_i \in W(k)$ relevant $y_i$. Posons :
$$F(X):=\prod_{i=1}^T(1-Y_i^pX)^{h_i}$$
Il est clair que $F$ rel\`eve $f$. V\'erifions que $F$ est de la forme
annonc\'ee. On a :
\begin{eqnarray*}
  F(X^p) & = & \prod_{i=1}^T(1-(Y_iX)^p)^{h_i} \\
   & = & \prod_{j=0}^{p-1} \prod_{i=1}^{T} (1-\zeta^kY_iX)^{h_i} \\
\end{eqnarray*}
%o\`u $\zeta$ est une racine $p$-i\`eme de l'unit\'e.

\noindent
Or, on a $(\prod_{i=1}^T(1-y_i x)^{h_i})=f(x)^{F^{-1}}$, donc
$\prod_{i=1}^T(1-y_i x)^{h_i}$ v\'erifie les hypoth\`eses du
lemme \ref{lemme : 2.12}
et il existe $\hat{Q}(X),\hat{R}(X),\hat{S}(X) \in W(k)[X]$ et $U\in
W(k)$ inversible tels que : 
 $$ \prod_{i=1}^{T} (1-Y_i(\zeta^jX))^{h_i}=(1+\zeta^jX\hat{Q}(\zeta^jX))^p+U(\zeta^jX)^m
 ( 1+\zeta^jX\hat{R}(\zeta^jX))+p\hat{S}(\zeta^jX)$$
ce que l'on peut \'ecrire aussi :
\begin{eqnarray*}
 \prod_{i=1}^{T} (1-Y_i(\zeta^jX))^{h_i} & = & (1+\zeta^jX\hat{Q}(\zeta^jX))^p(1+U(\zeta^jX)^m
 ( 1+\zeta^jX\tilde{R}(\zeta^jX)) \\
& + & p\tilde{S}(\zeta^jX))
\end{eqnarray*}
avec $\tilde{R}(X),\tilde{S}(X) \in W(k)[[X]]$. Ce qui donne :
\begin{eqnarray*}
F(X^p) & = & \prod_{j=0}^{p-1} (1+\zeta^jX\hat{Q}(\zeta^jX))^p \\
 &  & \prod_{j=0}^{p-1}(1+U(\zeta^jX)^m ( 1+\zeta^jX\tilde{R}(\zeta^jX)) p\tilde{S}(\zeta^jX))
\end{eqnarray*}   
que nous regardons modulo $p^2$. On a :
$$\prod_{j=0}^{p-1} (1+\zeta^jX\hat{Q}(\zeta^jX))^p \in
1+X^pW(k)[X^p]$$
et 
\begin{eqnarray*}
& & \prod_{j=0}^{p-1}(1+U(\zeta^jX)^m ( 1+\zeta^jX\tilde{R}(\zeta^jX))
p\tilde{S}(\zeta^jX)) \\
& = &  \prod_{j=0}^{p-1}(1+U(\zeta^jX)^m
( 1+\zeta^jX\tilde{R}(\zeta^jX)))+p\sum_{j=0}^p(\tilde{S}(\zeta^jX) \\
 & & \prod_{\substack{k \in \{0,\cdots,p-1\}\\k\neq
     j}}(1+U(\zeta^kX)^m ( 1+\zeta^kX\tilde{R}(\zeta^kX)))
 \;\mathrm{mod} \;p^2   
\end{eqnarray*}
Remarquons que la derni\`ere somme appartient \`a $(\zeta
-1)W(k)[[\zeta,X]] \cap W(k)[[X]]=pW(k)[[X]]$, donc :
$$p\sum_{j=0}^p(\tilde{S}(\zeta^jX) \prod_{\substack{k \in \{0,\cdots,p-1\}\\k\neq
     j}}(1+U(\zeta^kX)^m
 ( 1+\zeta^kX\tilde{R}(\zeta^kX)))=0\;\mathrm{mod} \;p^2  $$
Enfin, on a :
$$\prod_{j=0}^{p-1}(1+U(\zeta^jX)^m( 1+\zeta^jX\tilde{R}(\zeta^jX))) 
\equiv (1+U^pX^{pm}(1+X^p\tilde{R}^p(X))) \; \mathrm{mod} \; (\zeta
-1)$$
ainsi que :
$$\prod_{j=0}^{p-1}(1+U(\zeta^jX)^m( 1+\zeta^jX\tilde{R}(\zeta^jX)))$$
$$\in W(k)[[X^p]] \cap
(1+X^mW(k)[[X]])=1+(X^p)^{([\frac{m}{p}]+1)}W(k)[[X^p]]$$
Donc $F(X)$ est de la forme annonc\'ee.

%\begin{flushright}
%$\square$
%\end{flushright}
 Nous allons montrer que l'\'equation $Y^p=F(X)$ d\'efinit une courbe
 ayant bonne r\'eduction sur R relativement \`a la valuation de Gauss en
$T:=\lambda^{\frac{-p}{m}}X$.

En effet, si on pose $Y=\lambda Z+1+XQ(X)$ et
$T:=\lambda^{\frac{-p}{m}}X$ alors l'\'equation $Y^p=F(X)$ donne en
r\'eduction :
$$z^p-z=u^pt^m$$
Encore une fois, on a l'\'egalit\'e des genres des fibres
g\'eom\'etriques et sp\'eciales, ce qui assure la bonne r\'eduction.

On obtient ainsi $n$ rev\^etements $Y_i^p=F_i(X)$ ($1\leq i \leq n$)
de $\pp_{\mathrm{K}}$ qui ont simultan\'ement bonne r\'eduction pour la m\^eme
valuation de Gauss (l'\'equation en r\'eduction est
$z_i^p-z_i=u_it^m$). On consid\`ere le produit fibr\'e de ces
rev\^etements, apr\`es normalisation il induit un rev\^etement
$\mathcal{C} \rightarrow \pp_{R}$ galoisien de groupe
$(\Z/p\Z)^n$. De plus, la fibre sp\'eciale $C_s$ est int\`egre car les $u_i$ sont
lin\'eairement ind\'ependants sur $\fp$ (cf. remarque 3 du paragraphe
1.3). Il reste \`a voir que ce rev\^etement a bonne r\'eduction sur R.

On \'ecrit $m+1=qp^{n-1}$, $q \in \N^*$. Le degr\'e de la diff\'erente
sp\'eciale du compositum des $n$ extensions $z_i^p-z_i=u_it^m$ est :

\begin{eqnarray*}
d_s & = & (m+1)(p-1)(1+p+\cdots p^{n-1}) \\
& = & qp^{n-1}(p-1)(1+p+\cdots p^{n-1})
\end{eqnarray*}    

Notons $d_{\eta}$ le degr\'e de la diff\'erente du
rev\^etement $C_{\eta} \rightarrow \pp_{K}$. Ce rev\^etement n'est
ramifi\'e qu'en les points qui sont des rel\`evements des p\^oles des
formes diff\'erentielles, i.e au plus $T=q(1+\cdots +p^{n-1})$
points (voir la d\'emonstration du lemme \ref{lemme  : 1.6}). Les groupes d'inertie \'etant cycliques d'ordre $p$, on
obtient :
$$d_{\eta} \leq p^{n-1}(q(1+p+\cdots +p^{n-1}))(p-1)=d_s$$
On obtient la bonne r\'eduction en appliquant le crit\`ere local de
bonne r\'eduction donn\'e dans \cite{gre1}.
\begin{flushright}
$\square$
\end{flushright}

\noindent
\underline{Remarque} : Si on regarde cette derni\`ere action en
r\'eduction modulo l'id\'eal maximal de R, on trouve un
$(\Z/p\Z)^n$-torseur au dessus de $k[[t]]$ donn\'e par les \'equations
:
$$\left\{
  \begin{array}{ccc}
    z_1^p-z_1 & = & u_1t^m \\
     & \vdots & \\
    z_n^p-z_n & = & u_nt^m 
  \end{array}
\right.
$$
o\`u les $u_i$ sont $\fp$-ind\'ependants, car attach\'es \`a un espace
$L_{m+1,n}$ (cf. remarque 3 du paragraphe 1.3).

\subsubsection{D\'eformation des $(\Z/p\Z)^2$-torseurs}

On s'int\'eresse ici \`a la d\'eformation de $(\Z/p\Z)^2$-torseurs au
dessus de Spec$k[[t]]$
avec une g\'eom\'etrie \'equidistante; ceci impose de ne consid\'erer
que des extensions de $k[[t]]$ pour lesquelles les sous-extensions
interm\'ediaires ont des conducteurs \'egaux.
 
Dans le cas $m+1=p$, on sait d'apr\`es \cite{gre2} Th\'eor\`eme
III.3.1 que la g\'eom\'etrie qui apparait est \'equidistante. Ainsi on
a le th\'eor\`eme suivant corollaire du th\'eor\`eme \ref{theo : 2.9} :
\begin{theos}\label{theo : 2.13}
Soit $\mathrm{G}=(\Z/p\Z)^2$, $p \geq 3$ et $\mathrm{R}$ un anneau de valuation
discr\`ete dominant l'anneau des vecteurs de Witt de $k$. Supposons
que $G$ est un groupe d'automorphismes de $k[[z]]$ et que chacune des sous-extensions de $
k[[z]]^{\mathrm{G}}$, d'ordre $p$ a un conducteur \'egal \`a
$p$. Alors, on ne peut pas relever G en un groupe d'automorphismes de R$[[Z]]$.
\end{theos}

\noindent
\underline{\textit{D\'emonstration}} : On rappelle le crit\`ere de
rel\`evement donn\'e par \cite{gre1} Th\'eor\`eme I.5.1. Pour qu'il y ait
rel\`evement, il faut et il suffit que :

``Etant donn\'e deux extensions interm\'ediaires de la forme
$k[[z]]^{\mathrm{G}_1}$, $k[[z]]^{\mathrm{G}_2}$, on puisse relever
chacune de ces extensions en
$\mathrm{R}[[Z]]^{\mathrm{G}_i}/\mathrm{R}[[Z]]^{\mathrm{G}}$
telles que ces deux derniers rev\^etements aient exactement $(p-1)$
points de branchement en commun.''

Supposons que le rel\`evement soit possible et traduisons ce qui doit
se passer au niveau de la fibre sp\'eciale $\ddd$. Les \'equations des
deux rev\^etements interm\'ediaires sont de la forme 
$Y_1^p=F_1(X)$ et $Y_2^p=F_2(X)$. On a $m<p$, donc pour chacune de ces
extensions, $\ddd$ est une
droite projective sur laquelle on a $m+1$ points
\'equidistants qui correspondent aux sp\'ecialisations des points
fixes par l'automorphisme d'ordre $p$ correspondant (cf. \cite{gre2} th\'eor\`eme III.3.1). On a  donc
des fonctions $f_1,f_2$ (qui sont les r\'eductions de $F_1$ et $F_2$)
de la forme : 
$$f_1=\prod_{i=0}^p(1-x_ix)^{h_i} \;\;\; h_0=0,\; h_i\neq 0,\; \sum h_i=0$$
$$f_2=\prod_{i=0}^p(1-x_ix)^{h'_i} \;\;\; h'_p=0,\; h'_i\neq 0,\; \sum
h'_i=0$$
telles que d$f_1$ et d$f_2$ aient les conditions sur leurs diviseurs
\'enonc\'ees pr\'ec\'edemment. Cette \'ecriture traduit d\'ej\`a le
fait que l'on a $(p-1)$ points de branchement en commun (leurs
sp\'ecialisations sont $x_1,\cdots,x_{p-1}$). 

Posons $\omega_1=\frac{\mathrm{d}f_1}{f_1}$ et
$\omega_2=\frac{\mathrm{d}f_2}{f_2}$. Toute autre extension
interm\'ediaire est donn\'ee par une \'equation de la forme
$\mathrm{Y}^p=f_1^{\epsilon_1}f_2^{\epsilon_2}$,
($(\epsilon_1,\epsilon_2) \in \fp^2 -\{(0,0)\}$), donc donne naissance \`a
une diff\'erentielle
$\omega_j=\epsilon_1\frac{\mathrm{d}f_1}{f_1}+\epsilon_2\frac{\mathrm{d}f_2}{f_2}$ qui a
aussi $p$ p\^oles distincts et un z\'ero d'ordre $m-1$ en $0$.

On voit donc que dans le cas $m+1=p$, la possibilit\'e de relever
l'action du groupe G implique l'existence d'espaces $L_{p,2}$, ce qui
est d\'ementi par le th\'eor\`eme \ref{theo :2.9}.

\begin{flushright}
$\square$
\end{flushright}

\noindent
\underline{\textit{Remarque 1}} : Dans \cite{ber}, Bertin donne des
obstructions au rel\`evement d'actions de groupe. Le th\'eor\`eme
\ref{theo : 2.13} donne
de nouvelles obstructions qui sont de nature diff\'erentielle.
 
\vspace{1em}
\noindent
\underline{\textit{Remarque 2}} : Le th\'eor\`eme \ref{theo : 2.13} utilise juste le
premier r\'esultat du th\'eor\`eme \ref{theo : 2.9}. Pour les cas $m+1=2p$ ou $3p$,
on ne peut pas \'enoncer un th\'eor\`eme analogue car on n'a pas
forc\'ement une g\'eom\'etriquement \'equidistante. On peut juste dire
dans ces cas-l\`a que si le rel\`evement est possible, il doit faire
apparaitre une g\'eom\'etrie plus complexe. 

\vspace{1em}

Enfin nous consid\'erons le cas o\`u $p=2$ avec cette fois-ci un
conducteur quelconque. On a alors le th\'eor\`eme suivant :

\begin{theos}\label{theo : 2.14}
On consid\`ere une action de G$=(\Z/2\Z)^2$ comme groupe
d'automorphismes de $k[[t]]$ dans laquelle chacune des sous-extensions
de $k[[t]]^{\mathrm{G}}$ d'ordre $2$ a m\^eme conducteur (on note
$m+1=2n$ ce conducteur). Alors on peut d\'eformer cette action en une
action de G sur R$[[T]]$, o\`u R$=W(k)[\lambda^{\frac{1}{2n-1}}]$.
\end{theos}
 
La d\'emonstration est d\^ue \`a Ito et suit des indications
de M.Matignon. Nous la redonnons avec quelques modifications.

On a tout d'abord besoin du lemme suivant :
\begin{lemme}\label{lemme : 2.15}
Soit $X_1,\cdots ,X_n$ $\in W(k)$ deux \`a deux distincts et $U \in W(k)^*$. Alors il
existe $X_{n+1}, \cdots ,X_{2n}$ (avec $X_i\neq X_j$ d\`es que $1 \leq
i < j \leq 2n$) et $Q(X),R(X) \in W(k)[X]$ tels que le polyn\^ome :
$$F(X):=\prod_{i=1}^{2n}(1-X_iX)=(Q(X))^2+UX^{2n-1}+2R(X)$$
et tels que le rev\^etement de Spec$(W(k)[X])$ donn\'e par $Y^2=F(X)$
ait bonne r\'eduction.
\end{lemme}

\noindent
\underline{D\'emonstration :} Notons $x_i$ la r\'eduction de $X_i$
modulo l'id\'eal maximal de $W(k)$. D'apr\`es le th\'eor\`eme
\ref{theo : 2.8},
on peut trouver $x_{n+1}, \cdots, x_{2n}$ tels que \nolinebreak:
$$f(x):=\prod_{i=1}^{2n}(1-x_ix)=(q(x))^2+ux^{2n-1}$$
(o\`u $u$ est la r\'eduction de $U$).
Choisissons des rel\`evements $X_i$ de $x_i$ et posons :
$$\tilde{F}(X)=\prod_{i=1}^{2n}(1-X_iX)$$

%On conserve les m\^eme notations : on
%prend R$=W(k)[\pi]$ (avec $\pi^m=\lambda$) ; on
%dispose donc d'une forme diff\'erentielle
%$\omega=\frac{\mathrm{d}f}{f}$. On va alors montrer que l'on 
%peut trouver un polyn\^ome $F \in $ $W(k)[X]$, relevant $f$, tel
%que l'\'equation $Y^p=F(X)$ d\'efinisse une courbe
%ayant bonne r\'eduction par rapport \`a la valuation de Gauss en
%$T=\lambda^{-\frac{p}{m}}X$. On aura alors une action
%de $(\Z /p\Z)$ sur le disque ouvert $p$-adique.

%On sait que :
%$$\frac{f'}{f}=\frac{ux^{m-1}}{\prod\limits_{i=0}^m(1+x_ix)}
%\;\;\; u \in k^*$$
%Soit $\tilde{F}$ un rel\`evement de $f$ que l'on choisit sous la forme
%factoris\'ee suivante \nolinebreak:
%$$\tilde{F}(X)=\prod_{i=0}^m(1+X_iX)^{h_i}$$
%o\`u $X_i$ rel\`eve $x_i$. 

%Vu l'expression de $\frac{f'}{f}$, on en d\'eduit que $f$ est de la
% forme :
%$$f(x)=(1+xa(x))^p+\frac{u}{m}x^m(1+xr(x))$$
%avec $a(x),r(x) \in k[x]$. 
\noindent
$\tilde{F}$ est aussi de la forme :
$$\tilde{F}(X)=Q^2(X)+2R(X)+UX^{2n-1}$$ 
avec $Q=1+a_1X+\cdots +a_{n}X^{n} \in W(k)[X]$, $R=b_1X+\cdots
+b_{2n-1}X^{2n-1} \in W(k)[X]$. Ecrivons
$\tilde{F}$ en fonction du param\`etre $T:=(-2) ^{-\frac{2}{2n-1}}X$ :
\begin{eqnarray*}
\tilde{F}(T) & = & Q^2((-2) ^{\frac{2}{2n-1}}T) +2(b_1(-2) ^{\frac{2}{2n-1}}T+
\cdots +b_{2n-1}(-2) ^2T^{2n-1})+\\
& & (-2)^2UT^{2n-1}
\end{eqnarray*}
Posons $Y=-2Z+Q$; si le coefficient $(b_1(-2)^{\frac{2}{2n-1}}T+
\cdots +b_{m}(-2) ^2T^{m})$ est nul modulo
$2$, alors on a en r\'eduction :
\begin{eqnarray} 
  \frac{((-2) Z+Q)^2-Q^2}{(-2) ^p} & =  & UT^m \nonumber \\
  Z^2-Z & = & UT^m \nonumber
\end{eqnarray}
Ceci est suffisant pour avoir la bonne r\'eduction. En effet, le
rev\^etement d'\'equation $Y^2=\tilde{F}(X)$ est ramifi\'e en $2n$ points
(nombre de racines de $\tilde{F}$), donc le genre de la fibre g\'en\'erique
est $\frac{(2n-2)(2-1)}{2}$ (formule d'Hurwitz), c'est-\`a-dire le
m\^eme que celui de la fibre sp\'eciale.

On va donc chercher \`a modifier $\tilde{F}$. On \'ecrit
$\tilde{F}(X)=\prod_{i=1}^{2n}(1-X_iX)$. Posons 
$$F(X)=\prod_{i=1}^{2n}\left(1-X_iX-2\epsilon _iX\right)$$
o\`u $\epsilon_i=0$ si $i\leq n$ et ($\epsilon _i, i>n$) sont des constantes \`a d\'eterminer pour avoir
bonne r\'eduction. 
%On a :
%$$\left(1+X_iX+2\epsilon _iX\right)\equiv
%(1+X_iX)^{h_i}+p\epsilon _i X(1+X_iX)^{h_i-1}\;\;\;\mathrm{mod}\;
%[p^2]$$
%donc 
\begin{eqnarray}
  F(X) & = & \prod_{i=1}^{2n}\left(1-X_iX\right)\left(1+2\sum_{i=1}^{2n-1}\frac{\epsilon
  _iX}{1-X_iX}\right) \;\mathrm{mod} \; [4] \nonumber \\
  & = & (Q^2+2R)\left(1+2\sum_{i=1}^{2n}\frac{\epsilon
  _iX}{1-X_iX}\right) +UX^{2n-1}\;\mathrm{mod} \; [4,X^{2n},2X^{2n-1}]
  \nonumber \\
  & = & Q^2+2\left(Q^2\sum_{i=1}^{2n}\frac{\epsilon
  _iX}{1-X_iX}+R\right)+UX^{2n-1}\;\mathrm{mod} \; [4,X^{2n},2X^{2n-1}]
  \nonumber \\
  & = & Q^2+2(Q^2X\sum_{i=1}^{2n}\epsilon _i(1+X_iX+\cdots
  +(X_iX)^{2n-2})+R)\nonumber \\
& + & UX^{2n-1}\;\mathrm{mod} \; [4,X^{2n},2X^{2n-1}]
  \nonumber 
\end{eqnarray}
On va donc s'arranger pour que le terme 
$$Q^2X\sum_{i=1}^{2n}\epsilon _i(1+X_iX+\cdots
  +(X_iX)^{2n-2})+R)$$ soit nul modulo $2$.
Remarquons tout d'abord que si $k\geq n$, alors les termes en $X^k$
(\'ecrits en fonction du param\`etre $T$) sont nuls modulo
$2$. Il suffit donc de voir que l'on peut choisir les $\epsilon
_i$ de telle fa\c con que les termes en $X^k$ ($1\leq k\leq n-1$) de
l'expression :
$$Q^2X\sum_{i=1}^{2n-1}\epsilon _i(1+X_iX+\cdots
  +(X_iX)^{2n-2})+R)$$
 soient nuls.
Soit $\alpha _k$ le $k$-i\`eme terme de la s\'erie de Taylor de
$(-RQ^{-2})$ (i.e $(-RQ^{-2})=\sum_{k \geq 1}\alpha _k X^k)$. Alors la
condition que l'on vient d'\'enoncer se ram\`ene au syst\`eme :
$$\sum_{i=n+1}^{2n}\epsilon _iX_i^k=-\alpha _k \;\;\;\mathrm{pour} \;\;\;0\leq k
\leq n-2$$
qui a des solutions puisque c'est un syst\`eme de Vandermonde avec des
\'equations en moins.

\begin{flushright}
$\square$
\end{flushright}
 
Revenons \`a la d\'emonstration du th\'eor\`eme. Consid\'erons une
$(\Z/2\Z)^2$-extension $k[[z]]/k[[t]]$ telle que les sous-extensions
interm\'ediaires $C_i$ aient le m\^eme conducteur $m+1=2n$. Apr\`es un
changement de param\`etre $t$, on peut supposer que $C_1$ et $C_2$
sont donn\'ees par les \'equations :
$$\left\{
  \begin{array}{ccc}
     C_1 : y_1^2+y_1 & = & \frac{u}{t^{2n-1}} \\
     C_2 : y_2^2+y_2 & = & \frac{p(t)}{t^{2n-1}}
  \end{array}
\right.$$
avec $u\in k^*$ et $p(t)=1+p_1t+\cdots+p_{2n-2}t^{2n-2}$. D'apr\`es
\cite{gre1} Th I.5.1, il faut pouvoir relever $C_1$ et $C_2$ de fa\c
con \`a ce que ces deux rev\^etements aient exactement $n$ points de
branchements en commun.

Posons $t'=t(p(t))^{\frac{-1}{2n-1}}$. Alors les deux extensions
interm\'ediaires sont donn\'ees par :
$$\left\{
  \begin{array}{ccc}
     C_1:y_1^2+y_1 & = & \frac{u}{t^{2n-1}} \\
     C_2:y_2^2+y_2 & = & \frac{1}{t^{'2n-1}}
  \end{array}
\right.$$
Soit $T$ un param\`etre du disque ouvert relevant $t$ et
$T':=T(p(T))^{\frac{-1}{2n-1}}$ un param\`etre relevant $t'$ ( et $P(T)$
est un rel\`evement de $p(t)$). Si on \'ecrit $T'=\tau (T)$, alors
$\tau$ d\'efinit un automorphisme du disque ouvert Spec$W(k)[[T]]$. Notons
$X=2^{\frac{2}{2n-1}}T^{-1}$. Alors $\tau$ induit un automorphisme sur
le disque ferm\'e Spec$W(k)\{\{X^{-1}\}\}$ (rappelons que que les
\'el\'ements de $W(k)\{\{X^{-1}\}\}$ sont les s\'eries formelles de la
forme $\sum_{ \nu \geq 0} a_{\nu}X^{-\nu}$ avec $\lim_{\nu \rightarrow
  \infty}a_{\nu} =0$). Ce qui donne
$\tau(X^{-1})=X^{-1}P(2^{\frac{2}{2n-1}}X^{-1})^{\frac{-1}{2n-1}}$ et
$\tau$ est l'identit\'e en r\'eduction. Soit
$\tilde{C}_2:Y_2^2=1+\frac{4}{T^{'2n-1}}$ un rel\`evement de $C_2$ que
l'on peut re\'ecrire en choisissant de nouveaux param\`etres :
$$(Y'_2)^2=1-(X')^{2n-1}=\prod_{i=1}^{2n}(1-X'_iX')$$
Les id\'eaux $(1-X'_iX')$ d\'efinissent des points distincts
dans Spec$W(k)\{\{X^{-1}\}\}$. Posons $(1-X_iX):=\tau ^{-1}
(1-X'_iX')$. On applique alors le lemme pr\'ec\'edent aux points $X_1
\cdots X_n$, ce qui permet d'obtenir un rev\^etement d'\'equation :
$$\tilde{C_1}:(Y'_1)^2=A(X)^2+2B(X)+UX^{2n-1}$$
qui a bonne r\'eduction et qui a $n$ points de branchement en commun
avec $\tilde{C}_2$. Le rel\`evement souhait\'e est alors donn\'e par
la normalisation de $\tilde{C}_1 \times_{W(k)[[X]]} \tilde{C}_2$.

%\begin{flushright}
%$\square$
%\end{flushright}

\noindent
\underline{Remarque} : Dans le cas $p>2$ une g\'en\'eralisation du
th\'eor\`eme \ref{theo : 2.14} est un probl\`eme ouvert. On s'aper\c coit d\'ej\`a
au vu du th\'eor\`eme \ref{theo : 2.9} que la condition $p/m+1$ n'est pas
suffisante : il faut en plus l'existence d'espaces $L_{m+1,2}$. 
%\newpage 
%Or, on a le lemme suivant :
%\begin{lemme}
%Soit $c \in k$ et $F \in k[x]$ un polyn\^ome de la forme $\alpha
%x+F_0(x^p)$. Alors, pour une valeur ad\'equate de la constante $u$, la
%forme diff\'erentielle $\omega =\frac{u\mathrm{d}x}{F(x)^{p-1}-c}$ est logarithmique.
%\end{lemme}
%\underline{\textit{D\'emonstration}} : On utilise l'op\'erateur de
%Cartier ($\mathcal{C}\omega=\omega$) :
%\begin{eqnarray}
%\frac{u}{F(x)^{p-1}-c} & = &
%\frac{u(F(x)^{p-1}-c)^{p-1}}{(F(x)^{p-1}-c)^p} \nonumber \\
%& = & \frac{u\sum\limits_{j=0}^{p-1}F(x)^{(p-1)j}c^{p-1-j}}{(F(x)^{p-1}-c)^p} \nonumber
%\end{eqnarray}
%Or, on peut d\'emontrer facilement par r\'ecurrence sur $k$ que :
%$$(F(x)^{(p-1)j})^{(k)}=(-1)^kj(j+1)\cdots (j+k-1) \alpha^kF(x)^{(p-1)j-k}$$
%Donc :
%\begin{eqnarray}
%(u(F(x)^{p-1}-c)^{p-1})^{(p-1)} & = &
%u(-1)^{p-1}(p-1)!\alpha^{p-1}c^{p-2} \nonumber \\
%& = & -u\alpha^{p-1}c^{p-2} \nonumber 
%$\end{eqnarray}
%D'o\`u :
%$$\mathcal{C}\omega=\frac{(u\alpha^{p-1}c^{p-2})^{\frac{1}{p}}}{F(x)^{p-1}-c}$$
%Il suffit alors de prendre $u=(u\alpha^{p-1}c^{p-2})^{\frac{1}{p}}$
%pour avoir le r\'esultat escompt\'e.

\newpage
\nocite{*}

\bibliographystyle{amsalpha}
\bibliography{difflog}

\begin{flushright}

Guillaume Pagot\\
Laboratoire de th\'eorie des nombres \\
et algorithmique arithm\'etique,\\ Universit\'e de Bordeaux
  I, \\ 351, cours de la lib\'eration, \\ 33405 Talence Cedex

email : pagot@math.u-bordeaux.fr
\end{flushright}
\end{document}